\documentstyle[12pt,amssymb,amscd]{amsart}  

\newtheorem{thm}{Theorem}[section]
\newtheorem{cor}[thm]{Corollary}
\newtheorem{prop}[thm]{Proposition}
\newtheorem{lemma}[thm]{Lemma}

\theoremstyle{remark}
\newtheorem{remark}[thm]{Remark}
\newtheorem{example}[thm]{Example}

\theoremstyle{definition}
\newtheorem{defn}[thm]{Definition}

\numberwithin{equation}{section}

\newcommand{\bbR}{{\Bbb R}}

\newcommand{\bbZ}{{\Bbb Z}}

\newcommand{\cF}{{\cal F}}

\newcommand{\cP}{{\cal P}}
\newcommand{\cL}{{\cal L}}

\newcommand{\cA}{{\cal A}}

\newcommand{\Sh}{\operatorname{Sh}}

\newcommand{\Ker}{\operatorname{Ker}}

\newcommand{\D}{\operatorname{D}}

\newcommand{\Hom}{\operatorname{Hom}}
\newcommand{\Ext}{\operatorname{Ext}}

\newcommand{\shHom}{\underline{\operatorname{Hom}}}

\newcommand{\Sym}{\operatorname{Sym}}

\newcommand{\Or}{{\underline{o}}}
\newcommand{\End}{{\operatorname{End}}}
\newcommand{\St}{\operatorname{St}}
\newcommand{\Link}{\operatorname{Link}}

\newcommand{\printname}[1]
  {\smash{\makebox[0pt]{\hspace{-2.0in}\raisebox{8pt}{\tiny #1}}}}

\thanks
{The authors were supported in part by the CRDF grant RM1-2405-MO-02. The second author
was partially supported by the NSA grant MDA904-01-1-0020.}

\begin{document}

\title{Hard Lefschetz theorem and Hodge-Riemann relations for intersection cohomology of
nonrational polytopes}

\author{Paul Bressler}
\address{Department of Mathematics, University of Arizona, Tucson, AZ, 85721, USA}
\email{bressler@@math.arizona.edu}
\author{Valery A.~Lunts}
\address{Department of Mathematics, Indiana University,
Bloomington, IN 47405, USA}
\email{vlunts@@indiana.edu}

\begin{abstract}
The Hard Lefschetz theorem for intersection cohomology of nonrational
polytopes was recently proved by K. Karu [Ka]. This theorem implies the conjecture
of R. Stanley on the unimodularity of the generalized $h$-vector. In this paper we strengthen
Karu's theorem by introducing a canonical bilinear form $(\cdot ,\cdot )_{\Phi}$
on the intersection cohomology
$IH(\Phi)$ of a complete fan $\Phi$ and proving the Hodge-Riemann bilinear relations
for $(\cdot ,\cdot )_{\Phi}$.
\end{abstract}

\maketitle

\section{Introduction}

For an $n$-dimensional convex polytope $Q$ R.~Stanley ([S])
 defined a set of integers
$h(Q)= (h_0(Q), h_1(Q),\dots, h_n(Q))$ - the ``generalized $h$-vector'' -
which are supposed to be the intersection cohomology Betti numbers of the
toric variety $X_Q$ corresponding to $Q$. In case $Q\subset\bbR^n$
is a rational polytope the variety $X_Q$ indeed exists, and it is known ([S])
that  $h_i(Q)=\dim IH^{2i}(X_Q)$. Thus, for a {\em rational polytope $Q$},
the integers $h_i(Q)$ satisfy
\begin{enumerate}
\item $h_i(Q)\geq 0$,
\item $h_i(Q)=h_{n-i}(Q)$ (Poincar\'e duality),
\item $h_0(Q)\leq h_1(Q)\leq ... \leq h_{[n/2]}(Q)$
(follows from the Hard Lefschetz theorem for projective
algebraic varieties).
\end{enumerate}

For an arbitrary convex polytope (more generally for an Eulerian poset)
Stanley proved ([S], Theorem 2.4) the property 2 above. He conjectured that 1 and 3
also hold without the rationality hypothesis.  This conjecture drew attention for the
last fifteen years and many attempts had been made to prove it. A correct proof was recently
given by K. Karu in [Ka].
The proof is based on the theory developed in [BrLu] and [BBFK], which, in particular,
associates to an arbitrary convex polytope $Q$ its {\em intersection cohomology space} $IH(Q)$
with a Lefschetz operator $l$. In case $Q$ is rational, $IH(Q)=IH(X_Q)$ and $l$ is the
multiplication by the Chern class of an ample line bundle on $X$; so $l$ satisfies the
Hard Lefschetz theorem on $IH(Q)$. The similar Hard Lefschetz (HL) property of $l$ was
conjectured for an arbitrary $Q$ in [BrLu], [BBFK]. It was shown in [BrLu] the (HL) property
implies in particular the equality $\dim IH^{2k}(Q)=h_k(Q)$, hence implies the
Stanley's conjecture. The (HL) theorem was proved in [Ka].

The construction of the intersection cohomology $IH(Q)$ is based on the study of the
{\em equivariant geometry} of the (nonexistent in general) toric variety $X_Q$ and was
essentially given in [BeLu]. Namely, one first defines the "equivariant" intersection
cohomology
of $Q$ and then the ordinary $IH(Q)$. This equivariant part of the theory is indispensable.

More precisely, we work with the complete fan $\Phi$ which is dual to the polytope $Q$
(and denote $IH(\Phi)=IH(Q)$, etc.).
We consider $\Phi$ as a (finite) topological space with subfans being the open subsets.
There are two natural sheaves of rings on $\Phi$: the constant sheaf $A_{\Phi}$ ($A$ is the
ring of global polynomial functions on the space of the fan) and the "structure sheaf"
$\cA _{\Phi}$, such that the stalk $\cA _{\Phi ,\sigma}$ at the cone $\sigma \in \Phi$
consists of polynomial functions on $\sigma $. (In case the space of the fan is the
Lie algebra of a torus, the (evenly graded) algebra $A$ is naturally isomorphic to the
cohomology ring of the classifying space of the torus and there is a close connection
between $A_{\Phi}$-modules and equivariant sheaves on the corresponding toric variety $X$.)

The theory of $A_{\Phi}$-modules and $\cA _{\Phi}$-modules was developed in [BrLu]
and partially in [BBFK]. This
includes Verdier duality, "equivariant perverse sheaves", decomposition theorem for
subdivisions, etc. In particular there exists a {\em minimal sheaf} $\cL _{\Phi}$
which plays the role of the equivariant intersection cohomology complex. The minimal sheaf
is characterized by the property that it is an indecomposable (locally free)
$\cA _{\Phi}$-module
with the stalk  at the origin $\cL _{\Phi ,\Or}=\bbR$. Then the "equivariant"
intersection cohomology of $\Phi$ is the graded $A$-module $\Gamma (\Phi ,\cL _{\Phi})$ and
$$IH(\Phi):=\Gamma (\cL _{\Phi})/A^+\Gamma (\cL _{\Phi}),$$
where $A^+\subset A$ is the maximal ideal. It was proved in [BrLu], [BBFK] that
$IH(\Phi)$ satisfies the Poincar\'e duality: $\dim IH^{n-k}(\Phi)=\dim IH^{n+k}(\Phi)$.

Note that $\Gamma (\cL _{\Phi})$ (and hence $IH(\Phi)$) is a module over $\cA (\Phi)$ --
the algebra of piecewise polynomial functions on $\Phi$. The following theorem
was formulated as the main conjecture in [BrLu] and proved in [Ka].

\begin{thm} (HL) Assume that $l$ is a strictly convex piecewise linear function on $\Phi$.
Then the map
$$l^k:IH^{n-k}(\Phi)\to IH^{n+k}(\Phi)$$
is an isomorphism for all $k\geq 1$.
\end{thm}

We show in this paper that by assuming the (HL) theorem for fans of dimension $\leq n-1$ we
obtain a {\em canonical} pairing
$$(\cdot ,\cdot )_{\Phi}:IH(\Phi)\times IH(\Phi)\to \bbR (2n)$$
for fans $\Phi$ of dimension $n$. Then we prove the following theorem which is the analogue
of the Hodge-Riemann bilinear relations in the (intersection) cohomology of
algebraic varieties.

\begin{thm} (HR) Let $\Phi$ and $l$ be as in the (HL) theorem. Consider the primitive subspace
$$Prim _lIH^{n-k}(\Phi)=\Ker \{l^{k+1}:IH^{n-k}(\Phi)\to
IH^{n+k+2}\}.$$ Then for $0\neq a \in Prim _lIH^{n-k}(\Phi)$ we
have $$(-1)^{\frac{n-k}{2}}(a,l^ka)_{\Phi}>0.$$
\end{thm}

This theorem was also essentially proved by Karu, except he did
not have a canonical pairing $(\cdot ,\cdot)_{\Phi}$ and had to
make choices at each step. This ambiguity makes the proof
unnecessarily heavy and hard to follow. Actually the results of
this paper imply that the pairing in [Ka] is {\em independent} of
the choices made and coincides with the canonical one. For
completeness we present the proofs of theorems (HL) and (HR),
following the main ideas of [Ka]. The proof proceeds by induction
on the dimension and a reduction to the simplicial case which is
known by McMullen [McM] (see also [Ti]).

We follow the idea in [BBFK] and try to develop all the notions and
results not only for complete fans, but also for {\em
quasi-convex} ones. (A fan $\Delta $ is quasi-convex if $\Gamma
(\Delta ,\cL _{\Delta})$ is a free $A$-module and hence it makes
sense to define the intersection cohomology $IH(\Delta)$ the same
way as for complete fans.) In particular, we define the canonical
pairing $(\cdot ,\cdot)$ for quasi-convex fans and then show that
this pairing is compatible with various natural operations on fans
such as subdivisions, embeddings of fans, etc.

Let us briefly describe the contents of the paper. In the second
section we recall and collect some elementary facts about fans and
polytopes which are used later. In the third section we recall the
theory of sheaves on a fan according to [BrLu] and [BBFK], and
formulate the (HL) and (HR) theorems. Section 4 is a review the
"smooth" case corresponding to a simple polytope or a simplicial
fan. Here we review Timorin's work in case of polytopes and
Brion's work in case of fans. We then relate the two pictures in a
natural way so that Timorin's Poincar\'e duality is identified with
Brion's. In section 5 we define the canonical pairing on the
intersection cohomology of quasi-convex fans and then in section 6
we show that for simplicial fans our pairing equals $n!$ times
Brion's (or Timorin's). This implies the (HL) and (HR) theorems in
the simplicial case. Section 7 discusses compatibility of the
canonical pairing with natural operations of fans. This makes the
theory flexible and easy to use. In section 8 we obtain some
immediate applications. Section 9 contains the proof of theorems
(HL) and (HR). In section 10 we prove the K\"unneth formula for the
intersection cohomology and derive the (HL) and (HR) theorems for
the product of fans (these results are used in section 9).

We thank Kalle Karu for his useful remarks on the first version of this paper.
Vladlen Timorin informed us that the definition of the polytope algebra $A(P)$ in
section 4 and Proposition 4.3 are due to A.V. Pukhlikov and A.G. Khovanskii. He also
claims that Proposition 4.2 was known to Minkowski.

\section{Preliminaries on fans and polytopes}

Let us fix some terminology. Consider a linear space $V\simeq \bbR ^n$. A (convex) {\it cone}
in $V$ is the intersection of a finite number of closed half-spaces $\{ L _i\geq 0\}$ where
$L_i\in W:=V^*$. A {\it face} of
a cone $\sigma$ is the intersection $\{ L=0\}\cap \sigma$ where $L$ is a linear function
on $V$ which is nonnegative on $\sigma$. A {\it fan} $\Phi$ in $V$ is a finite collection
of cones with the following properties.

\noindent 1. If $\sigma \in \Phi$ and $\tau \subset \sigma$ is a face, then $\tau \in \Phi$.

\noindent 2. The intersection of two cones in $\Phi$ is a face of each.

\noindent 3. The origin $\Or \in \Phi$.

We denote by $\Phi ^{\leq k}\subset \Phi$ the subfan consisting of all cones of dimension
$\leq k$.

Given a fan $\Phi$ denote by $\vert \Phi \vert \subset V$ its {\it support}, i.e. the union
of all cones in $\Phi$. We call $\Phi$ {\it complete } if $\vert \Phi \vert =V$.

\begin{example} A fan $\Phi$ is {\it simplicial} if every cone $\sigma \in \Phi$ of dimension
 $k$ is a convex hull of $k$ rays.
\end{example}

For a subset $S\subset \Phi$ denote by $[S]$ the minimal subfan of $\Phi$ which contains $S$.
Let $\sigma \in \Phi$. We denote
$$\begin{array}{l}
\partial \sigma =[\sigma ]-\sigma,\\
\sigma ^0=\sigma -\vert \partial \sigma \vert,\\
\St (\sigma )=\{ \tau \in \Phi \vert \sigma \subset \tau \},\\
\partial \St (\sigma)=[\St (\sigma )]-\St (\sigma ),\\
\Link (\sigma)=\{ \tau \in \partial \St (\sigma )\vert \tau \cap \sigma = \Or \}\\
\end{array}
$$

Sometimes we will be more specific and write $\St _{\Phi}(\sigma )$ instead of $\St (\sigma )$.

We denote by $\langle \sigma \rangle$ the linear subspace of $V$ spanned by
a cone $\sigma $.

Let $\Phi$ and $\Phi ^\prime $ be fans in spaces $V$ and $V^\prime
$ respectively. Let $\gamma :\vert \Phi \vert \to \vert \Phi
^\prime \vert$ be a homeomorphism which is linear on each cone in
$\Phi$ and induces a bijection between $\Phi $ and $\Phi ^\prime$.
We call $\gamma $ a {\it pl-isomorphism} between $\Phi $ and $\Phi
^\prime$ and say that these fans are pl-isomorphic.

Let $V_1,$ $V_2$ be vector spaces with fans $\Delta $ and $\Sigma
$ respectively. Consider the product $V=V_1 \times V_2$ with the
fan $$\Phi =\{ \sigma +\tau \vert \sigma \in \Delta ,\tau \in
\Sigma \}. $$ We call $\Phi =\Delta \times \Sigma $ the product fan.

 Consider $A=\Sym V^*$ as a graded algebra,
where linear functions have degree 2.
 We denote by $\cA (\Phi)$ the (also evenly graded) algebra of piecewise polynomial functions
  on $\vert \Phi \vert$, i.e. $f\in \cA (\Phi)$ if $f_{\sigma }=f\vert _{\sigma}$ is a
   polynomial for each $\sigma \in \Phi$.

A piecewise linear function $l\in \cA ^2(\Phi)$ is called
 {\it strictly convex} if for any two cones $\sigma $ and $\tau $  of dimension $n$
  we have $l_{\sigma }(v)<l_{\tau}(v)$ for any $v\in \tau ^0$. A complete fan $\Phi$
  is called {\it projective}
  if it possesses a strictly convex piecewise linear function. Note that if $l$ is
  strictly convex, then so is $l+l^\prime$ for any $\l^\prime \in A^2$.

\begin{example} Let $P\subset W$ be a convex polytope of dimension $n$.
 For each face $F\subset P$ consider the cone $\sigma _F\subset V$ which is the dual to the cone in $W$
  generated by vectors $f-p$ with $f\in F$ and $p\in P$. The collection of cones $\sigma _F$
   is a complete fan in $V$, which is called the {\it outer normal fan } of $P$. We denote
    this fan $\Phi _P$. The assignment $F\mapsto \sigma _F$ is a bijective, order-reversing
     correspondence between faces of $P$ and cones of $\Phi _P$. Denote by $H_P$ the following
      function on $V$
$$H_P(v)=\max_{y\in P}\langle y,v\rangle.$$
Then $H_P$ is a piecewise linear function on $\Phi _P$, called the {\it support function}
 of $P$. It is strictly convex. Thus $\Phi _P$ is projective. Vice versa, given a
 projective fan  $\Phi$ with a strictly convex piecewise linear function $l$
 there exists a convex polytope $P\subset W$ such that $(\Phi ,l)=(\Phi _P ,H_P)$.
\end{example}

\begin{example} Consider the product fan $\Phi =\Delta \times
\Sigma$ as above. If $l_1$ and $l_2$ are piecewise linear strictly
convex functions on $\Delta $ and $\Sigma $ respectively, then
$l=l_1+l_2$ is strictly convex on $\Phi$.
\end{example}

\begin{remark} Our main results imply that a  strictly convex piecewise
linear function $l$ on a complete fan
behaves like the first Chern class of an ample line bundle on a projective
variety. This may cause some confusion, because of the two opposite notions
of strict convexity.
Namely, if $l=H_P$ is the support function of a
convex polytope on $\Phi =\Phi _P$, then $l$ is strictly {\it lower convex}:
 given two adjacent $d$-dimensional cones $\sigma ,\tau\in \Phi$, we have
$l_{\sigma}\vert _{\tau}<l_{\tau}$. In [Da] however, a piecewise linear
function which is the first Chern class of an ample line bundle, is
strictly {\it upper convex}. The reason for this discrepancy is the fact that
$\Phi _P$ is the {\it outer} normal fan of $P$, whereas Danilov uses the
{\it inner } normal fan.
\end{remark}

\medskip

Let $\Phi$ and $\Psi$ be fans. We call $\Psi $ a {\it subdivision} of $\Phi$ if
$\vert \Psi \vert =\vert \Phi \vert$ and every cone $\Psi$ is contained in a  cone
 in $\Phi$. In this case we define a map $\pi :\Psi \to \Phi$ so that
$\pi (\sigma )$ is the smallest cone in $\Phi$ which contains $\sigma $. Often we
 will refer to the map $\pi $ as a subdivision.

\begin{example} Given a fan $\Phi$ and a cone $\sigma \in \Phi$ of dimension $>1$
choose a ray $\rho \subset
\sigma ^0$. Define a new fan
$$\Psi :=\{\Phi -\St (\sigma)\}\cup \{ \tau +\rho \vert \tau \in \partial \St (\sigma )\}.$$
Then $\Psi $ is a subdivision of $\Phi$, called the {\it star} subdivision
(defined by $\rho $). Thus to obtain $\Psi $ from $\Phi$ one replaces $\St _{\Phi}(\sigma )$ by
$\St _{\Psi}(\rho)$.
\end{example}

\begin{remark} Given a star subdivision as above there is a one dimensional vector space of
 piecewise
linear functions on $\Psi$ supported on $\St (\rho)$. Namely, such a function is
determined by its restriction to $\rho $.
\end{remark}

\begin{defn} An edge (i.e. a 1-dimensional face) $\rho $ of a cone $\sigma $ is called
{\it free} in $\sigma$ if
all other edges of $\sigma $ are contained in one facet $\tau $, i.e.
$$\sigma =\tau +\rho.$$
A cone $\sigma $ is called {\it deficient} if it has no free edges.
\end{defn}

Note that a fan is simplicial if and only if it contains no deficient cones.

\begin{remark} In case of a star subdivision as in the previous example the $1$-dimensional
cone $\rho $ is free in every cone containing it.
\end{remark}

\begin{defn} For a fan $\Phi$ we define its "singular" subfan $\Phi ^s$ to be the minimal
subfan
which contains
all deficient cones in $\Phi$. So $\Phi ^s=\emptyset$ if and only if $\Phi$ is
simplicial.
\end{defn}

In the notation of Example 2.5 let $\sigma \in \Phi ^s$ be a cone of maximal dimension.
By Remark 2.8 the resulting fan $\Psi$ will contain a smaller number of deficient cones than
$\Phi$. This suggests a "desingularization" process of a fan, which ends when one obtains a
simplicial fan. This strategy is used in [Ka] and we also follow it.

\begin{defn} Let $\Phi $ be a fan and $\sigma \in \Phi$. We say that $\Phi$ has a local
product structure at $\sigma $ if
$$[\St(\sigma )]=\{ \xi =\tau _1 +\tau _2 \vert \tau _1 \in [\sigma], \tau _2 \in
\Link (\sigma) \}.$$
\end{defn}

\begin{example} Assume that a complete fan $\Phi $ has a local product structure
at $\sigma \in \Phi$.  Consider the projection $p: V\to
\overline{V}=V/\langle \sigma \rangle$. Then $\Psi=p(\Link (\sigma ))$ is a
complete fan in $\overline{V}$ and $[\St(\sigma)]$ is
pl-isomorphic to the product fan $\Psi \times [\sigma ]$.
\end{example}

\begin{lemma} Let $\Phi$ be a complete fan and $\sigma \in \Phi
^s$ - a maximal cone in $\Phi ^s$. Then $\Link (\sigma)$ is a
simplicial fan and $\Phi$ has a local product structure at $\sigma
$.
\end{lemma}

\begin{pf} {\it Claim:} Every cone $\tau \in \St(\sigma )$ is of
the form
 $$\tau =\sigma +\rho _1+...+\rho _{\dim (\tau)-\dim
(\sigma)},$$ where the edges $\rho _i$ are free in $\tau$.

Indeed, let $\tau \in \St (\sigma )$ be a cone of the minimal
dimension such that  $\tau = \sigma +\rho _1+...+\rho _s$ and
$\dim (\tau )<\dim (\sigma )+s$. By definition of $\sigma $ the
cone $\tau $ contains a free edge $\rho$, so that $\tau =\rho
+\xi$ for a cone $\xi$. If $\rho \subset \sigma$, then $\sigma
=\rho +\sigma \cap \xi$ and $\rho $ is free in $\sigma$ -- a
contradiction. So $\rho$ is one of $\rho _i$'s, say $\rho =\rho
_s$. But then $\sigma +\rho _1 +...+\rho _{s-1}$ is also a cone in
$\St(\sigma )$ and its dimension is less than $\dim (\sigma
)+(s-1)$ -- a contradiction with the minimality of $\tau$.  This
proves the first part of the claim. The freeness of each edge
$\rho _i$ is clear.

The lemma follows from this claim. Indeed, by definition every cone
in $\Link (\sigma )$ is a face of a cone in $\St (\sigma )$. The
claim implies that all cones in $\Link (\sigma )$ are simplicial.
This proves the first assertion of the lemma. The claim also
implies that for any face $\delta $ of $\sigma $ and any subset
$S$ of the cones $\{ \rho _i \}$ the cone $\delta +S$ is a face of
$\tau$. The local product structure follows.
\end{pf}

\begin{lemma} Let
$\Phi $ and $\sigma \in \Phi ^s$ be as in the previous lemma.
Choose a ray $\rho \in \sigma ^0$ and let $\Psi \to \Phi$ be the
corresponding star subdivision. If $\Phi$ is projective, so is
$\Psi$.
\end{lemma}

\begin{pf} Let $\dim (\sigma )=k$. Assume that $l\in \cA ^2(\Phi)$ is strictly convex on
$\Phi$.  Then $l$ is also in
 $\cA ^2(\Psi)$, but is not strictly convex on $\Psi$. Let $\tilde{l}$ be a piecewise
linear function on $\Psi$ supported on $\St (\rho)$, such that $\tilde{l}\vert _{\rho}<0$.
Consider
$$\hat{l}:=l+\epsilon \tilde{l}, \quad \text{for}\ \ 0<\epsilon <<1.$$
We claim that $\hat{l}$ is strictly convex on $\Psi$. Indeed,
let $\delta$ and $\theta $ be two cones
 of dimension $n$ in $\Psi$. If one of them is not in $\St (\rho)$, then
 $\hat{l}_{\delta}\vert _{\theta}<\hat{l}_{\theta}$, because the same is true for $l$
(and $\epsilon <<1$). So we may assume that $\delta, \theta \in
\St (\sigma )$ and, moreover, that they are contained in the same
cone $\tau \in \Phi$. The Claim in the proof of last lemma implies
that $$\tau =\sigma +\sum_{i=1}^{n-k}\rho _i,\quad \text{where}\ \
\rho _i\in \Link (\sigma ).$$ Therefore $$\delta =\delta \cap
\sigma +\sum_{i=1}^{n-k}\rho _i,\quad \theta =\theta \cap \sigma +
\sum_{i=1}^{n-k}\rho _i,$$ and $\dim (\delta \cap \sigma)=\dim
(\theta \cap \sigma )=k$. Now convexity of the boundary $\partial
\sigma $ implies that $\tilde{l}_{\delta}\vert
_{\theta}<\tilde{l}_{\theta}$ and hence  $\hat{l}_{\delta}\vert
_{\theta}<\hat{l}_{\theta}$, which proves the lemma.
\end{pf}

\begin{defn} A fan $\Delta $ is {\it purely $n$-dimensional} if every cone in
$\Delta $ is contained in a cone of dimension $n$. In this case we denote by
$\partial \Delta \subset \Delta $ the subfan whose support is the boundary of
the support of $\Delta $. Put $\Delta ^0=\Delta -\partial \Delta$.
\end{defn}

Unless stated otherwise all fans are assumed to be in $V$.

\section{Review of some constructions and results from [BrLu] and [BBFK]}

\subsection{}
We briefly recall the notions and results from our paper [BrLu]
which are relevant to this work. All $A$-modules are graded and a
morphism between $A$-modules is always homogeneous of degree $0$.
Let $A^+$ be the maximal ideal of $A$. For an $A$-module $M$ we
denote by $\overline{M}= M/A^+M$ the corresponding (graded)
$\bbR$-vector space. The shifted module $M(i)$ is defined by
$$M(i)^{k}=M^{i+k}.$$ Similarly we will consider shifts
$H^\bullet(i)$ of graded vector spaces $H^\bullet$.

The tensor product $\otimes $ always means $\otimes _{\bbR}$.

Let $\Phi $ be a fan in $V$. We consider the (finite) set $\Phi $ as a topological space
with open subsets being subfans. (In particular, $\Or \in \Phi$ is the unique open point).
Denote by $\Sh(\Phi)$ the category of sheaves of $\bbR$-vector spaces on $\Phi$.

There are two natural sheaves of rings on $\Phi $: the constant sheaf $A_{\Phi}$ and the
"structure sheaf" $\cA _{\Phi}$. The stalk $\cA _{\sigma }$ consists of polynomial functions
on $\sigma $. These sheaves are evenly graded, $\Gamma (\Phi ,A_{\Phi })=A$, $\Gamma
(\Phi ,\cA _{\Phi })=\cA (\Phi)$. Thus the canonical surjection of sheaves
$$A_{\Phi}\to \cA _{\Phi}$$
does not induce a surjection of global sections in general.

The category of $\cA _{\Phi}$-modules contains an important subcategory
$\frak{M}=\frak{M}(\Phi)$, consisting of flabby {\it locally free} sheaves of finite
type. Namely, an
$\cA _{\Phi}$-module $F$ belongs to $\frak{M}$ if it is a flabby sheaf and $F_{\sigma}$ is a
finitely generated free $\cA_{\Phi ,\sigma}$-module for each $\sigma \in \Phi$.
The category $\frak{M}$ is of finite type:
every object is a direct some of indecomposable ones. There exists a distinguished
indecomposable object $\cL =\cL _{\Phi}\in \frak{M}$ characterized by the property
$\cL _{\Or}=\bbR$. We call $\cL $ the {\it minimal} sheaf. It is an analogue of the
equivariant intersection cohomology sheaf on the toric variety $X_{\Phi}$
 corresponding to $\Phi$
(this variety exists only if the fan is rational). It is easy to
see that $\cL $ can be characterized as  an object in $\frak {M}$
such that $\cL _{\Or}=\bbR $ and the map $$\overline{ \cL _{\sigma
}}\to \overline{\Gamma (\partial \sigma ,\cL )}$$ is an
isomorphism for each cone $\sigma $.

The following (easy) result is the combinatorial analogue of the (equivariant)
decomposition theorem for the direct image of perverse sheaves under proper maps.

\begin{thm} (Decomposition theorem) Let $\pi :\Psi \to \Phi$ be a subdivision.
The direct image functor $\pi_*:\cA_{\Psi }-mod\to \cA _{\Phi}-mod$ takes the
category $\frak{M}(\Psi)$ to $\frak{M}(\Phi )$. In particular the sheaf $\pi_*\cL_{\Psi}$
contains $\cL_{\Phi}$ as a (noncanonical) direct summand.
\end{thm}

The last assertion of the theorem follows from the first one since
$(\pi_*\cL _{\Psi})_0=\cL_{\Psi ,0}=\bbR$.

\begin{defn} In the notation of the last theorem a map $\alpha :\cL _{\Phi}\to \pi_*\cL _{\Psi}$
is called {\it admissible} if
$\alpha _0=id:\bbR \to \bbR$.
\end{defn}

Since $\cL$ is flabby,
$$H^i(\Phi ,\cL)=0,\quad \text{for $i>0$}.$$

\begin{defn} [BBFK] A fan $\Phi$ is called {\it quasi-convex} if $H^0(\Phi ,\cL)$ is a
free $A$-module.
\end{defn}

A complete fan is quasi-convex. More precisely we have the
following theorem.

\begin{thm}[BBFK]A quasi-convex fan is purely $n$-dimensional. A purely
$n$-dimensional fan $\Phi$ is quasi-convex if and only if $\vert \partial \Phi \vert $
is a real homology manifold. In particular, $\Phi $ is quasi-convex, if $\vert \Phi \vert$
is convex.
\end{thm}

\begin{example} Let $\Phi$ be a complete fan and $\sigma \in \Phi
$. Then the subfan $[\St (\sigma)]$ is quasi-convex.
\end{example}

For a quasi-convex fan $\Phi$ we define its intersection cohomology space as
$$IH^\bullet (\Phi):=\overline{\Gamma (\Phi ,\cL)}.$$
Denote by $ih_i(\Phi )=\dim IH^i(\Phi)$ the $i$-th Betti number of $\Phi$.
If $\Phi$ is rational, then $IH(\Phi)=IH(X_{\Phi})$.

\subsection{Hard Lefschetz theorem}
Fix a projective fan $\Phi$ with a strictly convex piecewise linear function $l$.
The $A$-module $\Gamma (\cL _\Phi)$ is in fact a $\cA(\Phi)$-module;
so $IH(\Phi)$ is also such.  In particular $l$ induces a degree 2 endomorphism
$$l:IH^\bullet(\Phi)\to IH^{\bullet +2}(\Phi).$$
It is known that $IH^i(\Phi)=0$ unless $i$ is even and $i\in [0,2n]$; $ih_0(\Phi)=1$,
and $ih_{n-k}(\Phi)=ih_{n+k}(\Phi)$.
The following "Hard Lefschetz" theorem was conjectured in [BrLu].
It was recently proved in [Ka].

\begin{thm} (HL) For each $k>0$ the map
$$ l^k:IH^{n-k}(\Phi)\to IH^{n+k}(\Phi)$$ is an isomorphism.
\end{thm}

\begin{cor} For a projective fan $\Phi$ we have
$$ih_0(\Phi)\leq ih_2(\Phi)\leq...\leq ih_{2[d/2]}(\Phi).$$
\end{cor}

It was shown
in [BrLu] that the (HL) theorem  implies the Stanley conjecture on the unimodality
of the generalized $h$-vector.

\subsection{Duality} Let $\omega = \omega_{A/\bbR}=A\cdot \wedge ^nV^*$ be the dualizing
module for $A$. It is a free $A$-module of rank 1 generated in degree $2n$. The
{\it dualizing} sheaf on $\Phi$ is
$$D_{\Phi}=\omega _{\Or }[n],$$
(which is a complex concentrated in degree $-n$ and supported at the origin $\Or$.)

The Verdier duality functor
$$D:D_c^b(A_{\Phi}-mod)^{op}\rightarrow D_c^b(A_{\Phi}-mod)$$
is defined by
$$D(F)=R\shHom _{A_{\Phi}}(F,D_{\Phi}).$$
It has the following properties.

\begin{thm}([BrLu]) Let $\Phi$ be a fan.

a) There is a natural isomorphism of functors $Id\to D^2$. So $D$ is an anti-involution of
$D^b_c(A_{\Phi}-mod)$.

b) $D(\frak{M}(\Phi))=\frak{M}(\Phi).$

c) $D\cL \simeq \cL$.
\end{thm}

The last assertion of the theorem implies the Poincar\'e duality for $IH(\Phi)$ for a complete
fan $\Phi$.

For a deeper study of our duality functor we will need the notion of a cellular complex which
we recall next.

\subsection{Cellular complex $C^\bullet (\cdot)$} Choose an orientation of each (nonzero)
cone in the fan $\Phi$. Then to each sheaf $F$ on $\Phi$ we can associate its
{\it cellular complex}
$$C^0(F)\stackrel{\partial }{\rightarrow}C^1(F)\stackrel{\partial }{\rightarrow} ...
\stackrel{\partial }{\rightarrow}C^n(F),$$
where
$$C^i(F)=\bigoplus_{\dim(\sigma)=n-i}F_{\sigma},$$
and the differential $\partial $ is the sum of the restriction maps $F_\sigma \to F_\tau$
with plus or minus sign depending on whether the orientations of $\sigma $ and $\tau $ agree
or not.
Sometimes we will be more specific and write $C^\bullet _{\Phi}(F)$ for $C^\bullet(F)$.
Note that $C^\bullet(\cdot)$ is an exact functor from sheaves on $\Phi$ to complexes.

In particular, we get a functor
$$C^\bullet(\cdot ):D^b(A_{\Phi}-mod)\rightarrow D^b(A-mod).$$
Note that $C^\bullet(D_{\Phi})=\omega$ (a complex concentrated in degree 0).

The following proposition will be used extensively in this work.

\begin{prop} [BrLu]  For any fan $\Phi$ there exists a natural isomorphism
of functors $D_c^b(A_{\Phi}-mod)^{op}\to D_c^b(A-mod)$
$$\bbR \Gamma (\Phi ,D(\cdot))=\bbR \Hom ^\bullet _{A_{\Phi}}(\cdot ,D_{\Phi})\to \bbR
\Hom _A^\bullet(C^\bullet(\cdot),\omega),$$
induced by the functor $C^\bullet (\cdot)$.
\end{prop}

This proposition shows that the Verdier duality can be considered as the Borel-Moore
duality (see [BrLu]). Namely, for $F\in D^b_c(A_{\Phi}-mod)$ the assignment
$$\sigma \mapsto C^\bullet (F_{[\sigma ]})$$
defines the co-sheaf (of complexes) of sections of $F$ "with compact supports".
(Here $F_{[\sigma]}$ denotes the extension by zero of the restriction of $F$ to
$[\sigma ]$). It follows from the above proposition that
$$\bbR \Gamma ([\sigma] ,D(F))= \bbR
\Hom _A^\bullet(C^\bullet(F_{[\sigma ]}),\omega),$$
where the right hand side is the analogue of the topological Borel-Moore duality.

\medskip

The cohomology groups $H^i(C^\bullet (F))$ tend to be related to
the cohomology groups $H^i(\Phi ,F)$. For example, if $\Phi$ is
complete and the orientations of all $n$-dimensional cones agree,
then the natural map $\Gamma (\Phi ,F)\to C^0(F)$ induces an
isomorphism $\Gamma (F)\simeq H^0(C^\bullet(F))$. Moreover, the
other cohomology groups also coincide, i.e. we have a canonical
functorial quasi-isomorphism of complexes (for a complete fan
$\Phi$) $$\bbR \Gamma (\Phi ,F)\simeq C^\bullet(F)$$ ([BrLu]).
This shows that the cellular complex makes sense once we choose a
global orientation on $V$. In fact for our purposes it will be
necessary to choose a volume form on $V$. So from now on we make
the following assumption.

\medskip

\noindent{\bf Assumption 1.} We assume that a nonzero element $\Omega =\Omega _V\in
\wedge ^nV^*$ has been fixed.

\medskip

This choice determines a trivialization $\omega =A(2n)$. Thus in particular we obtain
an isomorphism $C^\bullet(D_{\Phi})=C^n(D_{\Phi})=A(2n)$. We want to see how this isomorphism
changes in a natural way once we make a different choice of the volume form. For this assume,
for simplicity, that $\Phi$ is complete.

Two different choices of orientations of cones in $\Phi$ produce two
functors $C_1^\bullet$ and $C_2^\bullet$ and there exists an isomorphism
$\phi :C_1^\bullet \to C_2^\bullet$. We want $\phi$ to be compatible with the
canonical isomorphisms $H^0(C_1^\bullet(F))=\Gamma (F)=H^0(C_2^\bullet(F))$.
This forces $\phi ^0:C^0_1\to C^0_2$ to be the identity map, and hence defines
uniquely the isomorphisms $\phi ^i:C^i_1\to C^i_2$. For example, if the two
choices of orientations are the same for all cones of dimension $<n$, and are
opposite on the maximal cones (i.e. only the global orientation is changed),
then $\phi ^i=-1$ for all $i> 0$. One can easily check that the sign
of $\phi ^n$ depends only on the comparison of the orientations of
maximal cones, i.e. $\phi ^n=1$ if the global orientations agree and
$\phi ^n=-1$ if they are opposite.

In particular, let $W$ be a vector space and $W_{\Or}$ be the sheaf on $\Phi$,
 which is equal to $W$ on $\Or$ and zero elsewhere. Choose a cellular
complex $C^\bullet$ and get an isomorphism $C^\bullet(W_{\Or})=C^n(W_{\Or})\simeq W[-n]$.
However, as was explained above  if we choose a different cellular complex
 with the opposite global orientation,
then the above isomorphism changes sign.

A certain ``cancellation'' of two ambiguities occurs in the previous
example if $W=\omega$. Indeed, our choice of the volume form $\Omega $ defines an
isomorphism $\omega =A(2n)$ and also an isomorphism $C^\bullet (D_{\Phi})=\omega$.
So we get an isomorphism
$$C^\bullet (D_{\Phi})=A(2n).$$
But the arguments above show that the volume form $-\Omega $ defines the same (!)
isomorphism. Hence this last isomorphism is independent of $\Omega $ up to multiplication
by a positive number.

\begin{remark} The above situation is similar to the fact that on a smooth
orientable connected manifold the top cohomology group with coefficients
in the orientation sheaf is canonically isomorphic to $\bbZ$.
\end{remark}

Let us apply Proposition 3.9 together with the last argument
to our main example of the sheaf $\cL$.

\begin{thm} Assume that the fan $\Phi$ is complete. Then the volume  form $\Omega$
defines a canonical
isomorphism of $A$-modules
$$\Gamma (D\cL)=\Hom _A(\Gamma (\cL ),A(2n)).$$
A different choice of a volume form defines the same isomorphism up to
multiplication by a positive real number.
\end{thm}

\begin{pf} Indeed, the complex $C^\bullet(\cL )$ is acyclic except in degree 0, and
$H^0(C^\bullet(\cL))=\Gamma (\cL)$ is a free $A$-module.
\end{pf}

\begin{remark} Note that both $\cL $ and $D\cL$ are $\cA _{\Phi}$-modules (not just
$A_{\Phi}$-modules). The $\cA (\Phi)$-module structure on $\Gamma
(D\cL)$ comes from such structure on $\Gamma (\cL)$ via the
isomorphism of the last theorem.
\end{remark}

\subsection{The Hodge-Riemann bilinear relations} In order to formulate the second main
theorem we need to make the following assumption.

\medskip

\noindent{\bf Assumption 2.} From now on we assume that the (HL)
theorem holds for all complete fans of dimension $<n$.

\medskip

\begin{lemma} Let $\Phi$ be a fan with the minimal
sheaf $\cL _{\Phi}$. Then for any cone $\tau \in \Phi$ of dimension $k\geq1$
the following hold:

1)The stalk $\cL _{\Phi ,\tau}$ is generated (as the $A$-module) in degrees
$<k$.

2)The costalk $\Gamma _{\{ \tau\}}\cL _{\Phi}:=\Ker(\cL_{\Phi,\tau}\to
\Gamma (\partial \tau ,\cL_{\Phi}))$ is generated in degrees $>k$.
\end{lemma}

\begin{pf} To simplify the notation assume that $k=n$. Choose a ray $\rho$
in the interior of $\tau$ and consider the projection $p:V\to \overline{V}
\simeq \bbR ^{n-1}$. Then $\overline{\partial\tau}:=
p(\partial \tau)$ is a complete fan in $\overline{V}$ and
$\partial \tau\subset V$ is a graph of a stictly convex
 piecewise linear function $\bar{l}$ on $\overline{\partial \tau}$. Let
$A_{\overline{V}}$ be the algebra of polynomial functions on $\overline{V}$.
Note then $A=A_{\overline{V}}[\bar{l}]$. The projection
$p:\partial \tau \to \overline{\partial\tau}$ defines an isomorphism of $A_{\overline{V}}$ -
modules
$$\Gamma (\partial \tau ,\cL_{\Phi})\simeq \Gamma (\overline{\partial \tau} ,
\cL_{\overline{\partial \tau}}).$$
So in particular $\Gamma (\partial \tau ,\cL_{\Phi})$ is a free $A_{\overline{V}}$-module.
 By our Assumption 2 above $\bar{l}$ acts on
$IH(\overline{\partial \tau})$ as a Lefschetz operator. This implies
both assertion of the lemma because the residue map
$$\overline{\cL_{\Phi,\tau}}\to
\overline{\Gamma (\partial \tau ,\cL_{\Phi})}$$
is an isomorphism.
\end{pf}

\begin{cor} $\End_{A_{\Phi}}(\cL_{\Phi})=\bbR$.
\end{cor}

\begin{pf} This follows from the previous lemma and from the fact that
$\cL _{\Phi ,\Or}=\bbR$ by induction on the dimension of cones in $\Phi$.
Indeed, let $\sigma \in \Phi$ be of dimension $d$
and consider the exact sequence of $A$-modules
$$0\to \Gamma _{\{ \sigma\}}\cL _{\Phi}\to \cL_{\Phi,\sigma}\to
\Gamma (\partial \sigma ,\cL_{\Phi})\to 0.$$
Since $\cL_{\Phi,\sigma}$ (resp. $\Gamma _{\{ \sigma\}}\cL _{\Phi}$) is generated in degrees
$<d$ (resp. $>d$) the identity endomorphism of $\cL _{[\partial \sigma ]}$ can be extended
to an endomorphism of $\cL _{[\sigma ]}$ in a unique way.
\end{pf}

\begin{lemma} There exists a canonical isomorphism of $\cA _{\Phi}$-modules
$$\epsilon _{\Phi}:\cL _{\Phi}\to D(\cL _{\Phi}).$$
\end{lemma}

\begin{pf} Since the sheaves $\cL _{\Phi}$ and $D(\cL _{\Phi})$ are
isomorphic, by the previous corollary it suffices to find a canonical
isomorphism of the stalks $\cL _{\Phi ,\Or}$ and $D(\cL _{\Phi})_{\Or}$.
But
$$D(\cL _{\Phi})_{\Or}=R\Hom _A(\cL _{\Phi,\Or},D_{\Phi ,\Or})=
R\Hom _A(\bbR ,\omega [n])=\Ext ^n_A(\bbR ,\omega [n]).$$
To compute the last Ext group take the canonical Koszul resolution
of the $A$-module $\bbR$
$$0\to \omega \to A\otimes \wedge ^{n-1}V^*\to ...\to A\otimes V^*\to A
\to \bbR.$$
Then $\Ext^n_A(\bbR ,\omega [n])=\Hom _A(\omega [n],\omega [n])=\bbR.$
\end{pf}

\begin{cor} Assume that the fan $\Phi$ is complete. The volume form $\Omega $
defines a canonical nondegenerate pairing
$$[\cdot ,\cdot]=[\cdot ,\cdot]_{\Phi}=
\Gamma (\cL _{\Phi})\times \Gamma (\cL _{\Phi})\to A(2d),$$
and hence a canonical nondegenerate pairing
$$(\cdot ,\cdot)=(\cdot ,\cdot)_{\Phi}:IH(\Phi)\times IH(\Phi)\to \bbR (2d).$$
A different volume form defines the same pairings up to multiplication by a positive real
number. These pairings are $\cA (\Phi)$-bilinear.
\end{cor}

\begin{pf} This follows from the Theorem 3.11, Remark 3.12 and Lemma 3.15.
\end{pf}

In section 5 below we will construct a similar canonical pairing
for quasi-convex fans.

The next theorem is the analogue of the Hodge-Riemann bilinear relations on the primitive
intersection cohomology of complex algebraic varieties.

\begin{thm} (HR) Let $\Phi$ be a projective fan with  a strictly convex piecewise
linear function $l$. For $k\geq 0$ put $$Prim_lIH^{n-k}(\Phi):=\{
b\in IH^{n-k}\vert l^{k+1}b=0\}.$$ Then the quadratic form
$Q_l(a)=(-1)^{\frac{n-k}{2}}(a,l^ka)_{\Phi}$ is positive definite
on $Prim_lIH^{n-k}(\Phi)$.
\end{thm}

\subsection{Intersection cohomology of pl-isomorphic fans and of $\St [\sigma]$}
Let $\Phi$, $\Psi $ be fans and $\gamma :\vert \Phi \vert \to
\vert \Psi \vert $ be a pl-isomorphism. Then there exists an
isomorphism of sheaves $\gamma ^*\cA _{\Psi}\simeq \cA _{\Phi}$.
Hence also an isomorphism of minimal sheaves $\gamma ^*\cL
_{\Psi}\simeq \cL _{\Phi}$ and of the global sections $\Gamma
(\Phi ,\cL _{\Phi})\simeq \Gamma (\Psi ,\cL _{\Psi})$. This last
isomorphism however does not preserve the structure of a module
over the global polynomial functions, unless $\gamma$ is induced by a
linear isomorphism of the ambient vector spaces. Thus, for
example, if the pl-isomorphic fans are quasi-convex, there is no
canonical isomorphism between $IH(\Phi)$ and $IH(\Psi)$. Note
however, that these graded vector spaces {\it are} isomorphic,
since the graded spaces $\Gamma (\cL _{\Phi})$ and $\Gamma (\cL
_{\Psi})$ have the same Hilbert function.

Assume that a complete fan $\Phi$ has a local product structure at
a cone $\sigma $. Let $p:V\to \overline{V}=V/\langle \sigma
\rangle$ be the projection, $\Phi _{\sigma}=p(\Link (\sigma ))$ -- the
complete fan in $\overline{V}$. Note that the projection $p$ makes
$\Gamma ([\St (\sigma )],\cL _{\Phi})$ an $\cA(\Phi _{\sigma})$-module.
Indeed, by the local product structure at $\sigma $ the projection
$p(\tau )$ of any cone $\tau \in [\St (\sigma )]$ is a cone in
$\Phi _{\sigma}$. In particular this makes $IH([\St (\sigma )])$ a $\cA
(\Phi _{\sigma})$-module.

\begin{lemma} There exist natural isomorphisms of $\cA
(\Phi _{\sigma})$-modules $$\Gamma ([\St (\sigma )],\cL
_{\Phi})\simeq \Gamma (\Phi _{\sigma},\cL _{\Phi}) \otimes \cL
_{\Phi ,\sigma},$$ $$IH([\St (\sigma )])\simeq IH(\Phi
_{\sigma})\otimes \overline{\cL _{\Phi ,\sigma}},$$ with the
trivial $\cA (\Phi _{\sigma})$-module structure on $\cL _{\Phi
,\sigma}$. In particular, if $\sigma $ is simplicial, i.e. $\cL
_{\Phi ,\sigma}=\cA _{\Phi ,\sigma }$, then $$IH([\St (\sigma )])=
IH(\Phi _{\rho}).$$
\end{lemma}

\begin{pf} Let $A_{\sigma }$ be the polynomial functions on
$\sigma $. Choose a projection $V\to \langle \sigma \rangle$ and
identify $A_{\sigma }$ as a subalgebra of $A$.

Let $p_{\sigma}: \St (\sigma )\to \Phi _{\sigma}$ denote the restriction of
the projection p to $\St (\sigma )$. The sheaf $A_{\sigma }\otimes
p^{-1}_{\sigma }\cL _{\Phi _{\sigma}}$ extended by zero to $\Phi$ satisfies
the definition of a {\it minimal sheaf based at} $\sigma $ in
[BrLu] (5.1). (It is denoted by $\cL _{\Phi }^{\sigma }$ in [BrLu].)
Note that the restriction $\cL _{\Phi}\vert _{\St (\sigma)}$
extended by zero to $\Phi$ is also a minimal sheaf based at
$\sigma$. By Proposition 5.2 in [BrLu] $\cL _{\Phi}\vert _{\St
(\sigma)}\simeq \cL ^{\sigma}_{\Phi}\otimes \overline{\cL
_{\Phi,\sigma}}$. Therefore $$\Gamma (\St (\sigma ),\cL
_{\Phi})\simeq \Gamma (\Phi _{\sigma},\cL _{\Phi _{\sigma}})\otimes A_{\sigma}\otimes
\overline{\cL _{\Phi ,\sigma}}=\Gamma (\Phi _{\sigma},\cL _{\Phi _{\sigma}})\otimes
\cL _{\Phi ,\sigma}.$$
This is an isomorphism of $\cA
(\Phi _{\sigma })$-modules. It remains to note that $$\Gamma ([\St (\sigma
)],\cL _{\Phi})=\Gamma (\St (\sigma ),\cL _{\Phi}).$$
\end{pf}

\section{Review of the ``smooth'' case: simple polytopes and
simplicial fans}

\subsection{Review of Timorin's work on simple polytopes}

The analogue of Hodge-Riemann bilinear relations for simple polytopes was proved
by McMullen in [McM]. Later Timorin gave a simpler proof in [Ti].
Here we recall his main results.

Consider the dual space  $W=V^*$.
Let $P\subset W$ be a convex polytope of dimension $n$. We assume that $P$ is
{\it simple}, i.e. at each vertex of $P$ there meet exactly $n$ faces of
dimension $n-1$.

\begin{defn} Two convex polytopes $P^\prime, P^{\prime\prime}
 \subset W$ are called
{\it analogous} if their hyperplane faces can be pairwise matched so that
the matched faces

1) have the same outward normal direction,

2) are analogous.

Any two segments on a line are analogous.
\end{defn}

Clearly being analogous is an equivalence relation. Let $\cP ^+(P)$
denote the collection of all polytopes in $W$ analogous to $P$.
If $P^\prime, P^{\prime\prime}\in \cP^+(P)$, then their Minkowski sum
$$P^\prime +P^{\prime\prime}=\{p^\prime+p^{\prime\prime}\vert
p^\prime\in P^\prime,p^{\prime\prime}\in P^{\prime\prime}\}$$
is also in $\cP^+(P)$. Also
$$\lambda P^\prime=\{ \lambda p^\prime\vert p^\prime\in P^\prime\},
\quad \lambda >0$$
is in $\cP^+(P)$. Thus $\cP^+(P)$ has a structure of a convex cone. It can
be complemented to a vector space $\cP(P)$ by considering formal
differences of polytopes.

 Choose a volume form $\Omega _W\in \wedge ^nV$ on $W$ to be the dual of $\Omega _V$
and  orient $W$ in such a way that
$$Vol(P):=\int\limits_P^{}\Omega _W>0.$$

\begin{prop} [Ti] The volume function $Vol$ defined on $\cP^+(P)$
extends to a homogeneous polynomial function of degree $n$ on
$\cP(P)$. We denote this extended function again by $Vol$.
\end{prop}

Let $P_1,...,P_m$ be the faces of $P$ of dimension $n-1$. Choose the
corresponding outward normal covectors $\xi _1,...,\xi _m\in V$. That is,
each $\xi _i\vert_P$ achieves its maximum exactly on $P_i$.
 Then we associate to $P$ the numbers
$H_1,...,H_m$ defined by
$$H_i:=\max_{p\in P}\{\xi _i(p)\}=\xi _i(P_i).$$

The functions $H_1,...,H_m$ form a system of linear
coordinates on $\cP(P)$.

\subsection{The polytope algebra $A(P)$} Let $Diff$ be the algebra of
differential operators with constant coefficients on $\cP(P)$. It is a
commutative polynomial ring with generators
$$\partial _i=\frac{\partial}{\partial H_i},\quad i=1,...,m.$$
Let $I\subset Diff$ be the ideal
$$I=\{d\in Diff\vert dVol=0\}.$$
Put $A(P):=Diff/I$. Since $Vol$ is a homogeneous polynomial, the ideal $I$
is graded and so is the ring $A(P)$:
$$A(P)=\oplus_{k=0}^{n}A_k(P).$$

The following proposition is almost obvious.

\begin{prop} [Ti] The formula $(\alpha ,\beta )_T=\alpha \beta
Vol$ defines a nondegenerate pairing
$$(\cdot,\cdot)_T:A_k(P)\otimes A_{n-k}(P)\to \bbR.$$
\end{prop}

\begin{cor} $\dim A_k=\dim A_{n-k}$.
\end{cor}

\begin{thm} [Ti] We have $\dim A_k=h_k(P)$ -- the $k$-th component of the
$h$-vector $h(P)$.
\end{thm}

\subsection{Presentation of the algebra $A(P)$}
We have a natural embedding of linear spaces
$$i:W\hookrightarrow Diff,\quad i(a)=\sum_{i=1}^m\xi_i(a)\partial _i.$$
Let $P_{i_1},...,P_{i_k}$ be different faces of $P$ of dimension $n-1$. Put
$$D_{{i_1}...{i_k}}:=\partial _{i_1}...\partial _{i_k}\in Diff.$$

\begin{thm} [Ti] 1) For all $a\in W$, $i(a)\in I$ (translation
invariance of the volume function).

2) If $P_{i_1}\cap ...\cap P_{i_k}=\emptyset$, then $D_{i_1...i_k}\in I$.

3) The ideal $I$ is generated by elements $i(a)$, $D_{i_1...i_k}$ as in
1), 2).
\end{thm}

\subsection{A basis for $A(P)$} Fix a general linear function $t$ on $W$.
For each vertex $p\in P$ let $index(p)$ be the number of edges
coming out of $p$ on which the maximum of $t$ is attained at $p$.
That is $index(p)$ is the number of edges which go down from $p$.
Consider the face $F(p)$ spanned by all these edges. We have $\dim
F(p)=index(p)$. Thus we associated one face to each vertex (the
interior of $P$ is associated to the ``highest'' vertex).

Note that each face $F\subset P$ of codimension $k$ is the intersection of
exactly $k$ hyperplane faces: $F=P_{i_1}\cap...\cap P_{i_k}$.
Put $D_F:=\partial _{i_1}...\partial _{i_k}\in Diff.$

\begin{thm} [Ti]
The monomials $\{D_{F(p)}\vert  \text{$p$ is a vertex in\ } P\}$
form a basis of $A(P)$. More precisely, monomials $\{D_{F(p)}\vert
index(p)=k\}$ form a basis of $A_{n-k}$.
\end{thm}

\subsection{The Lefschetz operator} Consider the operator
$$L=L_P=\sum_{i=1}^mH_i(P)\partial _i\in A_1(P).$$

\begin{thm} [Ti] The operator (of multiplication by)
$L$ acts as a Lefschetz operator on $A(P)$, i.e. for all $i\leq
n/2$ the map $$L^{n-2i}:A_i(P)\rightarrow A_{n-i}(P)$$ is an
isomorphism.
\end{thm}

\begin{thm} [Ti] Let $$PrimA_i(P)=\{a\in A_i(P)\vert L^{n-2i+1}a=0\}$$
be the primitive part of $A_i(P)$. Then the symmetric bilinear
form $\langle a , b \rangle  _i=(-1)^i(a,L^{n-2i}b)_T$ is positive
definite on $PrimA_i(P)$.
\end{thm}

\begin{thm} [Ti] We have $L_P^nVol=n!Vol(P)$, i.e. $(L_P^n,1)_T=
n!Vol(P)$.
\end{thm}

\subsection{The dual picture in terms of fans after Brion [Bri]}
Let $P\subset W$ be a simple polytope and $\Phi =\Phi _P\subset V$
be its outer normal fan (Example 2.2). Since $\Phi$ is simplicial,
$\cL_{\Phi}=\cA _{\Phi}$ and we denote
$$H^\bullet(\Phi)=IH^\bullet(\Phi).$$ Then $H^\bullet(\Phi)$ is a
(evenly) graded algebra
$$H^\bullet(\Phi)=\oplus_{k=0}^nH^{2k}(\Phi).$$

\begin{prop} [Bri] $\dim H^{2k}(\Phi)=h_k(P)$.
\end{prop}

\subsection{A pairing  $(\cdot ,\cdot)_B$ on $H^\bullet(\Phi)$}

Brion defines a (degree zero homogeneous) map of $A$-modules
$\zeta :\cA (\Phi)\to A(2n)$.
It follows from the last proposition that such
a map is unique up to a scalar factor.
Let us
recall the construction. For each $n$-dimensional cone $\sigma$,
denote by $\cF _{\sigma }$ the product of equations of the facets of
$\sigma $. Then $\cF _{\sigma}\in A$ is uniquely defined up to scalar
multiplication. One normalizes $\cF _{\sigma }$ as follows: the equations
of the facets are nonnegative on $\sigma $ and their wedge product is equal
$\pm \Omega _V$. Denote by $\phi _{\sigma }\in \cA (\Phi)$ the function such
that
$$\phi _{\sigma }(v)=\begin{cases}\cF_{\sigma }(v) & \text{if $v\in \sigma$}\\
0, & \text{otherwise}
\end{cases}
$$ Thus $\phi _{\sigma }\in \cA _{2n}(\Phi)$ and it vanishes
outside $\sigma $.

For $f\in \cA (\Phi)$ define $$\zeta(f):=\sum_{\dim(\sigma
)=n}\frac{f_{\sigma}}{\cF _{\sigma}}.$$

\begin{thm} [Bri] The map $\zeta$ is a well defined map $\zeta:
\cA (\Phi)\to A(2n)$
such that

1) $\zeta $ is $A$-linear;

2) $\zeta (\phi _{\sigma })=1$ for all $n$-dimensional $\sigma \in
\Phi$.
\end{thm}

Clearly, $\zeta (f)=0$ if $deg(f)<2n$.

The map $\zeta$ induces a symmetric pairing
$$[\cdot,\cdot]_B:\cA(\Phi)\times \cA(\Phi)\to A(2n),\quad
[a,b]_B=\zeta(ab).$$

We also obtain a nonzero linear function $$\bar\zeta
:H^{2n}(\Phi)\to A/A^+=\bbR(2n),$$ and hence a pairing $$(\cdot
,\cdot)_B:H^\bullet(\Phi)\times H^\bullet(\Phi)\to \bbR(2n),\quad
(x,y)_B=\bar\zeta(xy).$$

\begin{prop} [Bri] The pairings
$[\cdot,\cdot]_B$ and $(\cdot ,\cdot)_B$ are
nondegenerate.
\end{prop}

\begin{remark} Notice that if the form $\Omega _V$ is changed by a factor of
$r\in \bbR$, then the map $\zeta$ and hence the pairings $[\cdot,\cdot]_B$,
$(\cdot ,\cdot)_B$ are changed by the factor $\vert r\vert ^{-1}$.
\end{remark}

Recall the support function $H_P\in \cA (\Phi)$ of $P$ (Example
2.2): $$H_P(v)=\max_{x\in P}\langle x,v\rangle.$$ Note, that if
$P$ contains the origin in $W$, then $H_P$ is nonnegative.

\begin{thm} [Bri] We have $\zeta (H_P^n)=n!Vol(P)$.
\end{thm}

\subsection{Relation between the pictures of Timorin and Brion}

\begin{thm} There exists a natural isomorphism of algebras
$$\beta :A(P)\rightarrow H(\Phi),\quad \beta:A_k(P)\stackrel{\sim}{\to}
H^{2k}(\Phi),$$
such that $\beta (L_P)=H_P$, and $(a,b)_T=(\beta (a),\beta (b))_B$.
\end{thm}

\begin{pf} First note that there is a natural homomorphism of algebras
$\tilde{\beta}:Diff\to \cA (\Phi)$. Indeed, vectors in $V$ which
lie in the 1-dimensional cones in $\Phi$ are naturally linear functions
on the space $\cP (P)$. Hence a differential operator of order 1 on $\cP(P)$
defines a linear function on each 1-dimensional cone of $\Phi$. However,
since $\Phi$ is simplicial, every such function extends uniquely to a
piecewise linear function on $\Phi$. This defines the homomorphism
$\tilde{\beta}:Diff\to \cA (\Phi)$.

Notice that the generators $D_{i_1,...,i_k}$ of the ideal
$I\subset Diff$  lie in the kernel of $\tilde{\beta}$. Indeed, let
$P_i$ be a face of $P$ of dimension $n-1$, $\xi_i\in V$ -- a
corresponding outward normal covector, $\rho _i\in \Phi$ -- the
corresponding 1-dimensional cone in $\Phi$ (so $\xi _i\in \rho
_i$), $\partial _i\in Diff$ -- the corresponding derivation. Then
the piecewise linear function $\tilde{\beta}(\partial _i)$ takes
value 1 on $\xi _i$ and is zero on all other 1-dimensional cones.
Thus $\tilde{\beta}(\partial _i)$ is nonzero only on the star of
$\rho _i$. Now, it is clear that if the intersection of
$n-1$-dimensional faces $P_{i_1},...,P_{i_k}$ is empty, then so is
the intersection  of the stars of the corresponding 1-dimensional
cones $\rho _{i_1},...,\rho _{i_k}$. Hence
$$\tilde{\beta}(D_{i_1,...i_k})=\tilde{\beta}(\partial
_{i_1}\cdot... \cdot \partial _{i_k})=0.$$

Notice also that $\tilde{\beta}$ maps the subspace $i(W)\subset Diff$ to
the space of linear functions on $V$.

It follows from Theorem 4.6 that $\tilde{\beta}$ descends to a
ring homomorphism $\beta :A(P)\to H(\Phi)$. This homomorphism is
surjective, since $\tilde{\beta}$ is surjective (functions
$\tilde{\beta}(\partial _i)$ as above generate $\cA(\Phi)$). Since
$\dim A(P)=\dim H(\Phi)$, $\beta$ is an isomorphism.

It is clear that $\beta (L_P)=H_P$ and the last assertion follows
from Theorems 4.10 and 4.15 above.
\end{pf}

\subsection{Hard Lefschetz and Hodge-Riemann bilinear form for
simplicial fans} Let $\Phi$ be a complete simplicial fan in $V$ with a strictly
convex piecewise linear function $l$. Then there exists
a convex polytope $P\subset W$ such that $\Phi=\Phi_P$ and $l=H_P$.

The following theorems are immediate corollaries of Theorem 4.8,
4.9 and 4.16 above.

\begin{thm} For $\Phi$ and $l$ as above the multiplication by $l$ acts as
a Lefschetz operator on $H(\Phi)$, i.e.
$$l^k:H^{n-k}(\Phi)\stackrel{\sim}{\rightarrow}H^{n+k}(\Phi)$$ for
all $k\geq 1$.
\end{thm}

\begin{thm} Let $\Phi$, $l$ be as in the last theorem. For $k\geq 1$ put
$$PrimH^{n-k}(\Phi):=\{ b\in H^{n-k}(\Phi)\vert l^{k+1}b=0\}.$$
Then the symmetric bilinear form $(-1)^{\frac{n-k}{2}}(a,b)_B$ is
positive definite on $PrimH^{n-k}(\Phi)$.
\end{thm}

\section{Duality for quasi-convex fans}

\subsection{} Following the idea in [BBFK] we develop the duality for quasi-convex fans.
Recall that a fan $\Delta$ in $V$ is {\it quasi-convex} if $\Gamma
(\Delta ,\cL _{\Delta})=\Gamma (\cL )$ is a free $A$-module.

Let $\Delta $ be a quasi-convex fan in $V$, $\Delta ^0:=\Delta -\partial \Delta$.
Recall that $\partial \Delta$ and $\Delta ^0$ are open and closed subsets of $\Delta $
respectively. For a sheaf $F$ on $\Delta $ we denote by $\Gamma _{\Delta ^0}F$ the
global sections of $F$ supported on $\Delta ^0$. That is
$$\Gamma _{\Delta ^0}F=\Ker \{ \Gamma (\Delta ,F)\to \Gamma ( \partial \Delta ,F)\}.$$

Consider the standard short exact sequence of sheaves $$0\to \cL
_{\partial \Delta} \to \cL \to \cL _{\Delta ^0}\to 0$$ and the
induced short exact sequence of cellular complexes $$0\to
C^\bullet (\cL _{\partial \Delta})\to C^\bullet (\cL )\to
C^\bullet (\cL _{\Delta ^0})\to 0.$$ It is proved in [BBFK] (Thms
4.3 and 4.11) that in the natural commutative diagram
$$\begin{array}{ccc} \Gamma _{\Delta ^0}\cL & \hookrightarrow &
\Gamma (\cL )\\ \downarrow & & \downarrow \\ H^0(C^\bullet (\cL ))
& \to & H^0(C^\bullet (\cL _{\Delta ^0}))
\end{array}
$$
the vertical maps are isomorphisms. Moreover, the complexes $C^\bullet (\cL)$ and
$C^\bullet (\cL _{\Delta ^0})$ are acyclic in degrees other than zero.

Recall that $IH^\bullet(\Delta ):=\overline{\Gamma (\cL )}$. One
also defines $$IH ^\bullet(\Delta ,\partial \Delta ):=
\overline{\Gamma _{\Delta ^0}\cL}.$$

It was shown in [BBFK] that there exists a (noncanonical) duality
between the free $A$-modules $\Gamma (\cL)$ and $\Gamma _{\Delta
^0}$ and hence a duality between the graded vector spaces
$IH(\Delta)$ and $IH(\Delta ,\partial \Delta)$. We are going to
make this duality canonical and to show that it is compatible with
various natural constructions, such as embeddings of fans,
subdivisions, etc.

\subsection{Definition of the pairing
$[\cdot ,\cdot]_\Delta :\Gamma (\cL)\times \Gamma _{\Delta ^0}\cL
\to A(2n)$}

By copying the proof of Theorem 3.11 and using the previous
remarks we obtain the following proposition.

\begin{prop} The functor $C^\bullet (\cdot)$ induces a natural isomorphism of $A$-modules
$$\Gamma (\Delta , D\cL )=\Hom _A(\Gamma _{\Delta ^0}\cL
,A(2n)).$$ This isomorphism is independent of the choice of the
volume form $\Omega$ up to multiplication by a positive real
number.
\end{prop}

Combining this isomorphism with the canonical isomorphism
$\epsilon _{\Delta} :\cL _\Delta \to D\cL _\Delta$ we obtain the
following analogue of Corollary 3.16.

\medskip

\begin{prop} There exists a natural isomorphism of $A$-modules
$$\Gamma (\cL) \to \Hom _A(\Gamma _{\Delta ^0},A(2n)),$$ hence
nondegenerate pairings $$[\cdot ,\cdot ]=[\cdot ,\cdot ]_\Delta
:\Gamma (\cL)\times \Gamma _{\Delta ^0}\cL \to A(2n)$$ and
$$(\cdot ,\cdot)=(\cdot ,\cdot )_\Delta : IH (\Delta ) \times
IH(\Delta ,\partial \Delta) \to \bbR (2n). $$ These pairings are
$\cA (\Delta)$-bilinear and are independent of the choice of the
volume form up to multiplication by a positive real number.
\end{prop}

\begin{cor} The graded vector spaces $IH(\Delta )$, $IH(\Delta ,\partial \Delta)$ are
concentrated in degrees from $0$ to $2n$. Also $IH^{2n}(\Delta
,\partial \Delta)\simeq \bbR(2n)$.
\end{cor}

\begin{pf} Since the sheaf $\cL$ is zero in negative degrees and $(\cL)^0=\bbR _{\Delta}$,
both assertions follows from the nondegeneracy  of the pairing $(\cdot ,\cdot)$.
\end{pf}

\begin{remark} In case the fan $\Delta $ is complete, i.e. $\Delta ^0=\Delta$, we have
$\Gamma _{\Delta ^0}(\cL )=\Gamma (\Delta ,\cL )$ and the pairing
of the last proposition coincides with the one in Corollary 3.16.
\end{remark}

\section{All pairings coincide in the simplicial case}
Let $\Delta $ be a simplicial quasi-convex fan in  $V$. Since $\cL
_{\Delta}=\cA _{\Delta}$ we may define a pairing on $\Delta $
using Brion's functional $\zeta $. Namely define $$[\cdot
,\cdot]_B:\Gamma (\Delta ,\cA _{\Delta })\times \Gamma _{\Delta
^0}\cA _{\Delta}\to A(2n), \quad [a,b]_B:=\zeta(ab).$$

\medskip

\begin{prop} [BBFK] The pairing $[\cdot ,\cdot ]_B$ is well defined and is perfect.
\end{prop}

\medskip

This also follows from the next proposition.

\medskip

\begin{prop} We have $[a,b]=n![a,b]_B$ for $a\in \Gamma (\Delta ,\cA _{\Delta})$,
$b\in \Gamma _{\Delta ^0}\cA _{\Delta}$.
\end{prop}

\begin{pf} Since the pairing $[\cdot ,\cdot ]$ is $\cA(\Delta)$-bilinear, $[a,b]=[1,ab]$.
Define
the linear functional
$$s:\Gamma _{\Delta ^0}\cA _{\Delta} \to A(2n),\quad s(a):=[1,a].$$
It suffices to prove that $s=n!\zeta$.
We will describe the map $s$ explicitly and then compare it with $\zeta$.

Recall
the Koszul resolution of the structure sheaf $\cA _{\Delta}$.
Namely, define the sheaf of vector
spaces $\Omega ^1=\Omega_{\Delta}^1$ as follows:
$$\Omega ^1_{\sigma}=\sigma ^{\perp}\subset V^*.$$
Put $\Omega ^i:=\wedge ^i\Omega ^1$. Thus the sheaf $\Omega ^i$ is supported
on the subfan $\Delta ^{\leq n-i}$; it is a sheaf of graded vector spaces which
are concentrated in degree $2i$.

 We have the canonical Koszul resolution of
the structure sheaf $\cA _{\Delta}$:
$$0\to A_{\Delta}\otimes \Omega ^n\to ...\to A_{\Delta}\otimes \Omega ^1\to
A_{\Delta}\to \cA _{\Delta}\to 0.$$
Note that $A_{\Delta}\otimes \Omega ^n=\D_{\Delta}$.

From the definition of the canonical morphism
$\epsilon _\Delta :\cA _{\Delta}\to D\cA _{\Delta}$ (Lemma 3.15)
it is clear that $\epsilon _{\Delta }$ is the projection of the complex
$$(*)\quad \quad 0\to A_{\Delta}\otimes \Omega ^n\rightarrow
...\rightarrow
A_{\Delta}\otimes \Omega ^1\rightarrow A_{\Delta}\to 0$$
on its leftmost nonzero term. Thus the map $s:\Gamma _{\Delta ^0}\cA _{\Delta}\to A(2n)$
 coincides with the isomorphism
$\Gamma _{\Delta ^0}\cA \to H^0(C^\bullet(A_{\Delta}\otimes \Omega ^\bullet))$
(given by the embedding $\Gamma _{\Delta ^0}\cA
\hookrightarrow C^\bullet(A_{\Delta}\otimes \Omega ^\bullet)$)  followed by the
projection $H^0(C^\bullet(A_{\Delta}\otimes \Omega ^\bullet))\to C^\bullet (D_{\Delta})=A(2n).$

First, recall the formula for the Koszul differential
$$\nu :A_{\Delta}\otimes \Omega ^j\to A_{\Delta}\otimes \Omega ^{j-1}:$$
$$\nu (g\otimes dx_1\wedge ...\wedge dx_j)= \sum_i(-1)^igx_i\otimes
dx_1\wedge ...
\wedge \hat{dx}_i\wedge ... \wedge dx_j.$$

The double complex $C^\bullet(A_{\Delta}\otimes \Omega ^\bullet)$ looks like
$$\begin{array}{ccccc}
...& & ...& & ...\\
\partial \uparrow & & \partial \uparrow & & \uparrow \\
C^2(A_{\Delta}\otimes \Omega ^2) & \stackrel{\nu }{\rightarrow} &
C^2(A_{\Delta}\otimes \Omega ^1) & \stackrel{\nu }{\rightarrow} &
C^2(A_{\Delta})\\
& & \partial \uparrow & & \partial \uparrow \\
& & C^1(A_{\Delta}\otimes \Omega ^1) & \stackrel{\nu }{\rightarrow} &
C^1(A_{\Delta})\\
& & & & \partial \uparrow \\
& & & & C^0(A_{\Delta})
\end{array}$$

Here $\nu $ is the Koszul differential and $\partial $ is the cellular
complex differential.
The rows are exact except at the rightmost column.
Choose $g\in \Gamma _{\Delta ^0}\cA _{\Delta}$. Then
$g$ defines a chain $\{ g_{\sigma }\}\in C^0(\cA _{\Delta })=C^0(A_{\Delta})$.
Put
$$\tilde{s} (g):=(-1)^{[n+1/2]}\nu ^{-1}\cdot \partial \cdot \nu ^{-1} \cdot
\partial \cdot ...
\cdot \nu ^{-1}\cdot \partial (\{g_{\sigma}\})\in
C^d(A_{\Delta}\otimes \Omega ^n)=A\cdot \Omega.$$
(The first, third, fifth, etc. composition $\nu ^{-1}\cdot \partial$ is taken
with the negative sign). It is clear that
$$\tilde{s}(g)=s(g)\Omega.$$

Let us understand the map $\tilde{s}$. Let $\xi \subset \tau$ be cones in
$\Delta$ of codimension $k$ and $k+1$ respectively. Define
$$\epsilon(\tau ,\xi)=\begin{cases}
1, & \text{if orientations of $\tau$ and $\xi$
agree};\\
-1, & \text{otherwise}.
\end{cases}$$
Then the map on stalks $\nu ^{-1}\cdot \partial :
(A_{\Delta}\otimes \Omega ^k)_{\tau}\rightarrow (A_{\Delta}\otimes \Omega ^{k+1})_{\xi}$
is equal to the wedging $\epsilon(\tau ,\xi)\frac{dy}{y}\wedge \cdot$
followed by restriction from $\tau $ to $\xi$,
where $y$ is any nonzero linear function on $\tau $, s.t. $y\vert_{\xi}=0$.

Choose a cone $\sigma \in \Delta$ of dimension $n$. Choose linear
functions $x_1,x_2,...,x_n\in A_{\sigma}$ such that

1) $x_i>0$ in the interior of $\sigma$;

2) $\prod x_i=0$ on the boundary $\partial \sigma$;

3) $dx_1\wedge ... \wedge dx_n=\Omega$. (This is possible assuming $n>1$. We omit here case
$n=1$, since it can be done directly.)

Let $f\in \Gamma _{\Delta ^0}\cA _{\Delta}$ be defined as follows
$$f_{\tau }=\begin{cases} \prod x_i, & \text{if $\tau =\sigma $}\\
                           0, & \text{ otherwise.}
\end{cases}
$$
In particular, $f=\phi _{\sigma }$ in the notation of section 4 above.

The following lemma implies the proposition.

\begin{lemma} We have $\tilde{s}(f)=n!\Omega$, or, equivalently,
$s(f)=n!$.
\end{lemma}

\noindent{\it Proof of lemma.} The element $\tilde{s}(f)\in C^n(A_{\Delta}\otimes
\Omega ^n)$ is a sum of $n!$ terms corresponding to the $n!$ choices of
a complete flag of faces of $\sigma$.

\noindent{\it Claim 1.} {\it The $n!$ summands in $\tilde{s}(f)$ are
equal. }

\noindent{\it Claim 2.} {\it Each summand is equal to $\Omega$. }

Clearly the claims imply the lemma. It is also clear that the only issue in
proving the claims is the sign. Indeed, by applying $n$ times the formula above for
$\nu ^{-1}\cdot \partial$
we find that each of the $n!$ terms is $\pm \Omega$.

Let us start with the flag $\sigma _n=\Or \subset \sigma _{n-1}\subset ...
\subset \sigma _0=\sigma$, where $\sigma _k=\{x_{n-k+1}=...=x_n=0\}
\cap \sigma$. Assume that $\sigma _k$ is oriented by the form
$\omega _k=dx_1\wedge ...\wedge dx_{n-k}$, ($\omega _n=1$).
Then by definition $\epsilon(\sigma _{k},
\sigma _{k+1})\omega _{k+1}=-dx_{n-k}\wedge \omega _{k+1}$.
Indeed, $-dx_{n-k}$
is the outward normal covector of $\sigma _{k+1}$ in $\sigma _k$. The
sign $\epsilon(\sigma _{k}, \sigma _{k+1})$ is negative for $k+1=n,n-2,...$.
Thus $\epsilon =-1$ for
$[n+1/2]$ pairs $\sigma _{k+1}\subset \sigma _k$. Hence the summand of
$\tilde{s} (f)$ corresponding to the above flag (with above orientaions) is equal
to $(-1)^{2[n+1/2]}\Omega=\Omega$.

It remains to proof that the summand is independent of the flag and of the
orientations of the cones of dimension $0<k<n$ in the flag.

First note that if the flag is the same, but the orientation of one of
the cones, say $\sigma _k$, is changed ($0<k<n$), then
$\epsilon(\sigma _{k-1}, \sigma _k)$ and $\epsilon(\sigma _k,\sigma _{k+1})$
change sign, so the corresponding summand of $\tilde{s} (f)$ remains unchanged.

Now assume that the flag is changed in one place: for some $0<k<n$ the cone
$\sigma _k$ is replaced with $\tau _k:=\{x_{n-k}=x_{n-k+2}=...=x_{n}=0\}$.
Let $\tau _k$ be oriented by the form $\omega _k ^\prime =dx_1\wedge ...
\wedge dx_{n-k-1}\wedge dx_{n-k+1}$. (Orientations of other cones
remain the same). Then
$$-dx_{n-k+1}\wedge \omega_{k+1}=(-1)^{n-k}\omega ^\prime _k,$$
$$-dx_{n-k}\wedge \omega_k^\prime=(-1)^{n-k}\omega _{k-1}.$$
Hence $\epsilon(\sigma _{k-1}, \tau _k)=\epsilon(\tau _k,\sigma _{k+1})$.
This changes the total number of minus signs among the $\epsilon$'s by $1$,
but
the composition
$$\nu ^{-1}\cdot \partial \cdot \nu ^{-1}\cdot \partial:(A_{\Delta}\otimes
\Omega ^{k-1})_{\sigma _{k-1}}\to (A_{\Delta}\otimes \Omega ^k)_{\tau _k}\to
(A_{\Delta}\otimes \Omega ^{k+1})_{\sigma _{k+1}}$$
involves wedging with $dx_{n-k+1}\wedge dx_{n-k}$ as opposed to the
previous  wedging with $dx_{n-k}\wedge dx_{n-k+1}$. Thus the contribution to
$\tilde{s} (f)$ does not change.
This proves the claims and the lemma.
\end{pf}

\begin{cor} The (HL) and (HR) Theorems holds for simplicial fans. In particular,
they holds if $n\leq 2$.
\end{cor}

\begin{pf} Indeed, this is an immediate consequence of Theorems 4.17, 4.18 and
Proposition 6.2.
\end{pf}

\section{Compatibility of the pairing with natural operations on fans}

\subsection{Compatibility of the pairing with embedding of fans}
Let $\Sigma $ be a fan in $V$ and $\Delta \subset \Sigma $ a subfan. Assume that both
$\Sigma $ and $\Delta $ are quasi-convex.
Note that $\cL _{\Sigma}\vert _{\Delta}=\cL _{\Delta}$,
hence we get a (surjective) restriction homomorphism
$$r=r_{\Sigma, \Delta}: \Gamma (\Sigma, \cL _{\Sigma})\to \Gamma (\Delta, \cL _{\Delta}).$$
Also note that the closed subset $\Delta ^0\subset \Delta$ remains closed
in $\Sigma$ and $\Delta ^0\subset \Sigma ^0$. We obtain the (injective) map
$$s=s_{\Delta ,\Sigma}:\Gamma _{\Delta ^0}\cL _{\Delta}=\Gamma _{\Delta ^0}
\cL _{\Sigma}\to \Gamma _{\Sigma ^0}\cL _{\Sigma}.$$

\medskip

\begin{prop} The maps $r$ and $s$ are adjoint. That is
$$[r(a),b]_{\Delta}=[a,s(b)]_{\Sigma}$$
for  $a\in \Gamma (\Sigma, \cL _{\Sigma})$, $b\in \Gamma _{\Delta ^0}\cL _{\Delta}$.
\end{prop}

\begin{pf} Note that the restriction $\vert _{\Delta}$ of sheaves from
$\Sigma $ to $\Delta $ commutes with the duality.
Moreover, the following diagram commutes
$$\begin{array}{ccc}
\cL _{\Sigma} & \stackrel{\epsilon _{\Sigma}}{\rightarrow} & D(\cL _{\Sigma})\\
\vert _{\Delta} \downarrow & & \downarrow \vert _{\Delta}\\
\cL _{\Delta} & \stackrel{\epsilon _{\Delta}}{\rightarrow} & D(\cL _{\Delta}).
\end{array}
$$

Also note that the restriction $\Gamma (\Sigma ,D\cL _{\Sigma})\to
\Gamma (\Delta ,D\cL _{\Delta})$ is induced by the embedding of cellular complexes
$C^\bullet(\cL_{\Delta})\hookrightarrow C^\bullet(\cL_{\Sigma })$, which in turn
induces the inclusion $s:\Gamma _{\Delta ^0}\cL _{\Delta}\to \Gamma _{\Sigma ^0}\cL _{\Sigma}.$
That is the following diagram commutes
$$\begin{array}{ccc}
 \Gamma (\Sigma ,D\cL _{\Sigma}) & = &
 \Hom _A(\Gamma _{\Sigma }\cL _{\Sigma},A(2d))\\
\vert _{\Delta} \downarrow  &  & \downarrow s^*\\
\Gamma (\Delta ,D\cL _{\Delta})) & =
& \Hom _A(\Gamma _{\Delta }\cL _{\Delta},A(2d)).
\end{array}
$$
This proves the proposition.
\end{pf}

\subsection{Compatibility of the pairing with subdivisions of fans}

Let $\pi :\Psi \to \Phi$ be a subdivision of fans. Choose an admissible morphism
$\alpha :\cL _{\Phi}\to \pi_*\cL _{\Psi}$ (3.2). Since
$\cL _{\Phi}$ is a direct summand of $\pi_*\cL _{\Psi}$ it follows by the rigidity of
$\cL _{\Phi}$ (3.14) that $\alpha $ is a split injection.

\medskip

\begin{thm} Let $\Delta $ be a quasi-convex fan, $\pi :\Theta \to \Delta $ - a subdivision
(hence $\Theta $ is also quasi-convex).
Choose an admissible morphism $\alpha :\cL _{\Delta}\to \pi_*\cL _{\Theta }$ (3.2). It induces
morphisms
$$\alpha : \Gamma _{\Delta ^0}\cL _{\Delta}\to \Gamma _{\Delta ^0}(\pi _*\cL _{\Theta})=
\Gamma _{\Theta ^0}\cL _{\Theta},$$
$$\alpha : \Gamma (\Delta ,\cL _{\Delta })\to \Gamma (\Delta ,\pi _*\cL _{\Theta})=
\Gamma (\Theta ,\cL _{\Theta }).$$
Then for  $a\in \Gamma (\Delta ,\cL _{\Delta })$, $b\in \Gamma _{\Delta ^0}\cL _{\Delta}$, :
$$(\alpha (a),\alpha (b))_{\Theta}=(a,b)_{\Delta}.$$
\end{thm}

\begin{pf} First we want to "descend" the duality from $\Theta $ to $\Delta $. It was proved
in [BrLu] that the duality functor commutes with the direct image functor $\bbR \pi _*$. Here
 we want to describe this commutation in a way that is compatible with Proposition 3.9.

 Let $F\in D^b_c(A_{\Theta }-mod)$ and consider the cellular complexes
 $C^\bullet _{\Theta }(F)$
 and $C^\bullet _{\Delta }(\pi_* F)$.
 Recall that the $n$-dimensional cones in $\Theta $ and $\Delta $ are oriented by $\Omega _V$.
 We can orient the other cones in $\Delta $ and $\Theta $ in a compatible way.
 Namely, if $\dim (\tau )=\dim(\pi (\tau ))$
 we want the orientations of $\tau $ and $\pi (\tau )$
 to be compatible. Otherwise choose an orientation of $\tau $ at random.
Then there is a natural morphism of complexes
$$\phi (F):C^\bullet _{\Delta }(\pi _*(F))\to C^\bullet _{\Theta }(F),$$
which is compatible with identifications
$$H^0(C^\bullet _{\Theta }(F))=\Gamma (\Theta , F)=\Gamma (\Delta ,\pi _*F)=
H^0(C^\bullet _{\Delta }(\pi _*F)).$$

The functorial map $\phi$ above allows us to define a morphism of functors
$\gamma : \bbR \pi _*\cdot D\to D\cdot \bbR \pi _*$ in the following way.
Fix $F\in D^b_c(A_{\Theta}-mod)$. Let $F\to J^\bullet$ be its (finite) injective resolution
and choose an injective resolution $\omega \to I^\bullet$ of the $A$-module $\omega $.
For a subfan $\Phi \subset \Delta$ we have
$$D\cdot \bbR \pi _*(F)(\Phi)=\Hom _A(C^\bullet _{\Delta}((\pi _*J^\bullet)_{\Phi}),I^\bullet).$$
On the other hand, for a subfan $\Psi \subset \Theta $
$$D(F)(\Psi)=\Hom _A (C^\bullet _{\Theta}(J^\bullet _{\Psi} ),I^\bullet).$$
Note that $D(F)$ is a complex of flabby sheaves, hence $\bbR \pi _*D(F)=\pi _*D(F)$.
We define
$$\gamma (F)(\Phi): \bbR \pi _*\cdot D(F)(\Phi )\to D\cdot \bbR \pi _*(F)(\Phi)$$
 to be
$$\phi ^*  (J^\bullet ):\Hom _A(C^\bullet_{\Theta}(J^\bullet _{\pi ^{-1}(\Phi)}),I^\bullet)\to
  \Hom _A(C^\bullet _{\Delta }((\pi _*J^\bullet )_{\Phi} ),I^\bullet).$$

\medskip

\begin{lemma} $\gamma $ is an isomorphism of functors.
\end{lemma}

\begin{pf} Fix a cone $\sigma \in \Delta$. It suffices to show that for an injective
$A_{\Theta }$-module $J$ the map of complexes
$$\phi (J):C^\bullet ((\pi _*J)_{[\sigma ]})\to C^\bullet (J_{\pi ^{-1}([\sigma ])})$$
is a quasi-isomorphism. For simplicity of notation we may assume that $\dim (\sigma )=n$.
Put $G=\pi _*J$. Consider the exact sequence of sheaves
$$0\to G_{\partial \sigma} \to G_{[\sigma ]}\to G_\sigma \to 0.$$
Applying Proposition 3.6 of [BrLu] to sheaves $G_\sigma $ and $G_{\partial \sigma }$
on fans $[\sigma ]$ and $\partial \sigma $ respectively we obtain natural quasi-isomorphisms
$$\bbR \Gamma ([\sigma ],G_\sigma)=C_{\Delta}^\bullet (G_{\sigma}), \quad
\bbR \Gamma (\partial \sigma ,G_{\partial \sigma })=C_{\Delta}^\bullet
(G_{\partial \sigma })[1].$$
The sheaf $G$ is injective, hence the exact triangle
$$C_{\Delta}^\bullet (G_{[\sigma ]})\to C_{\Delta}^\bullet (G_\sigma )\to
C_{\Delta}^\bullet (G_{\partial\sigma })[1]$$
is quasi-isomorphic to the (middle part of the) short exact sequence
$$(1)\quad \quad \quad 0\to \Gamma _\sigma G \to \Gamma ([\sigma ],G)
\to \Gamma (\partial \sigma ,G) \to 0.$$

The same arguments show that the triangle
$$C_{\Theta}^\bullet (J_{\pi ^{-1}([\sigma ])})\to
C_{\Theta}^\bullet (J_{\pi ^{-1}(\sigma )})\to
C_{\Theta}^\bullet (J_{\pi ^{-1}(\partial \sigma )})[1]$$
is quasi-isomorphic to the (middle part of the) short exact sequence
$$(2)\quad \quad 0\to \Gamma _{\pi ^{-1}(\sigma)} J \to \Gamma (\pi ^{-1}([\sigma ]),J)
\to \Gamma (\pi ^{-1}(\partial \sigma ) ,J) \to 0.$$
The short exact sequences (1) and (2) are isomorphic and the isomorphism $\Gamma _{\sigma}G
\simeq
\Gamma _{\pi ^{-1}(\sigma )}J$ coincides with the map
$\phi (J)$ under the above quasi-isomorphisms.
\end{pf}

\medskip

\begin{remark} As in the proof of Lemma 3.15 it is easy to see that the stalks
$(\pi _*D \cL _{\Theta})_{\Or}$ and $(D \pi _*\cL _{\Theta })_{\Or}$
are canonically isomorphic to $\bbR$. Under this identifications $\gamma (\cL _{\Theta})_{\Or}=
id$.
\end{remark}

\medskip

 Consider the diagram of sheaves
$$\begin{array}{ccccc}
\pi _*\cL _{\Theta} & \stackrel{\pi _*\epsilon _{\Theta}}{\to } &
\pi _* D\cL _{\Theta} & \stackrel{\gamma}{\to } & D\pi _*\cL _{\Theta} \\
\alpha \uparrow & & & & \downarrow D(\alpha )\\
\cL _{\Delta } & & \stackrel{\epsilon _{\Delta}}{\to } & & D\cL _{\Delta}.
\end{array}
$$

\medskip

\begin{lemma} The above diagram commutes.
\end{lemma}

\begin{pf} Note that the stalks at $\Or$ of all the sheaves are canonically isomorphic to
$\bbR$ and all the morphisms are equal to $id$ at $\Or$. Hence the diagram commutes by
the rigidity of $\cL _{\Delta}$.
\end{pf}

It follows from the last two lemmas that the following diagram commutes
$$\begin{array}{ccccc}
\Gamma (\Delta ,\pi _*\cL _{\Theta}) & \stackrel{\pi _*\epsilon _{\Theta}}{\rightarrow} &
\Gamma (\Delta ,\pi_* D\cL_{\Theta}) & = &
\bbR \Hom _A(C^\bullet_{\Theta}(\cL _{\Theta}),A(2d))\\
 & & \gamma \downarrow & & \downarrow \phi (\cL _{\Theta})^*\\
 \alpha \uparrow & & \Gamma (\Delta ,D\pi_*\cL _{\Theta}) & = &
 \bbR \Hom _A(C^\bullet_{\Delta}(\pi_*\cL_{\Theta}),A(2d))\\
 & & D(\alpha)\downarrow & & \downarrow C^\bullet_{\Delta}(\alpha ) ^*\\
 \Gamma (\Delta ,\cL _{\Delta}) & \stackrel{\epsilon _{\Delta}}{\rightarrow} &
 \Gamma (\Delta ,D\cL _{\Delta}) & = & \bbR \Hom _A(C^\bullet _{\Delta }(\cL _{\Delta}),A(2d)).
 \end{array}
 $$

 Finally note the commutativity of the natural diagram
 $$\begin{array}{ccc}
 \Gamma _{\Delta ^0}\cL _{\Delta} & \stackrel{\alpha }{\rightarrow} & \Gamma _{\Theta ^0}
 \cL _{\Theta } \\
 \downarrow & & \downarrow \\
 C^\bullet_{\Delta}(\cL _{\Delta}) & \stackrel{\phi (\cL _{\Theta })\cdot C^\bullet _{\Delta}
 (\alpha)}{\rightarrow} & C^\bullet_{\Theta}(\cL _{\Theta}).
 \end{array}
 $$
 The theorem follows.
\end{pf}

\begin{cor} For a complete fan $\Phi$ the pairing
$$[\cdot ,\cdot]:\cA(\Phi)\times \cA(\Phi)\to A(2d)$$
is symmetric.
\end{cor}

\begin{pf}  This is true if the fan $\Phi$ is simplicial (Proposition 6.2). For a
general fan $\Phi$ take a subdivision $\pi:\Psi \to \Phi$ where $\Psi $ is
simplicial and use the last theorem.
\end{pf}

\begin{cor} Let $\Phi$ be a complete fan. Let $a,b\in \Gamma (\cL_{\Phi})$
have disjoint supports. Then $[a,b]_{\Phi}=0$.
\end{cor}

\begin{pf} Same as that of the last corollary.
\end{pf}

\subsection{Local-global compatibility of the pairing}

Assume that $\Delta =\St (\rho )$ for a 1-dimensional cone $\rho
\in \Delta ^0$. Denote by $p:V\to \overline{V}:=V/\langle \rho
\rangle$ the projection along $\rho$. Then $\Phi :=p(\partial
\Delta)$ is a complete fan in $\overline{V}$. We want to relate
the pairings on $\Delta $ and on $\Phi $.

Denote $B=\Sym \overline{V}^*$. By Lemma 3.18 there is a canonical isomorphism of
$A$-modules $A\otimes _B\Gamma (\Phi _\rho ,\cL _{\Phi})\simeq
\Gamma (\Delta ,\cL _{\Delta})$, hence an
identification $IH(\Phi  )=IH(\Delta)$.

Let $x_n$ be a linear function on $V$ which is positive on the interior $\rho ^0$ of
$\rho$. We have $A=B[x_n]$. The pairing on $IH(\Phi )$ is defined once we choose a volume
form $\Omega _{\overline{V}}$ on $\overline{V}$. Choose it to satisfy
$$\Omega _V=dx_n\wedge \Omega _{\overline{V}}.$$

Let $\psi =x_n-f=0$ be the equation of the boundary $\partial \Delta$, where $f$ is
a piecewise linear function on $\Phi $. Multiplication by $\psi$ maps $\Gamma (\Delta ,
\cL _{\Delta})$ to $\Gamma _{\Delta ^0}\cL _{\Delta }$, hence it maps $IH(\Delta )$ to
$IH(\Delta ,\partial \Delta)$.

\medskip

\begin{prop} Let $a,b \in IH(\Phi)=IH(\Delta)$. Then
$$n(a,b)_{\Phi}=(a,\psi(b))_{\Delta}.$$
\end{prop}

\begin{pf} Step 1: Reduction to the simplicial case. Let $\overline{\pi }:\Psi \to \Phi$
be a subdivision, such that the fan $\Psi $ is simplicial. It induces the subdivision
$\pi :\Theta \to \Delta$, such that $\Theta =\St (\rho)$ and $p(\partial \Theta)=\Psi$.
Then $\Theta $ is also simplicial. Choose an admissible embedding
$\overline{\alpha }:\cL _{\Phi}\to \overline{\pi}_*\cL _{\Psi}$. It induces an admissible
embedding $\alpha :\cL _{\Delta}\to \pi _*\cL _{\Theta}$. By Theorem 7.2 the induces maps on
global sections are isometries (and $\alpha$ commutes with multiplication by $\psi$). Thus
it suffices to prove the proposition for $\Theta$, i.e. we may assume that $\Delta $ is
simplicial.

Step 2. By Proposition 6.2 above for $y\in \cA (\Delta)$,
$z\in \Gamma _{\Delta ^0}\cA _{\Delta}$
$$[y,z]_{\Delta }=n!\zeta (yz),$$
where $\zeta :  \Gamma _{\Delta ^0}\cA _{\Delta} \to A(2n)$ is the Brion
functional determined by the volume form
$\Omega _V$. Similarly, for $s,t\in \cA (\Phi )$
$$[s,t]_{\Phi}=(n-1)!\bar{\zeta} (st),$$
for the functional $\bar{\zeta} :\cA(\Phi)\to B(2(n-1))$, determined by
$\Omega _{\overline{V}}$. Thus
we must check that for a function $g\in \cA(\Phi)$ we have
$$\bar{\zeta} (g)=\zeta(\psi g).$$
It suffices to check the last equality for $g=\phi _{\sigma}$ (section 4),
where $\sigma \in \Phi$ is of dimension $n-1$. Then $\psi g=\phi _{\tau }$,
where $\tau \in Star ^0(\rho)$ is the preimage of $\sigma$ (by our choice of the volume
forms). So $\bar{\zeta}(g)=1=\zeta(\psi g)$ by Theorem 4.12.
\end{pf}

\medskip

\begin{cor} Multiplication by $\psi$ induces isomorphisms $\Gamma (\cL _{\Delta})\to
\Gamma _{\Delta ^0}\cL _{\Delta}$, $IH(\Delta)\to IH(\Delta ,\partial\Delta)$.
\end{cor}

\begin{pf} This follows from the last proposition and the nondegeneracy of the pairings
$(\cdot ,\cdot)_{\Delta}$ and $(\cdot ,\cdot)_{\Phi}$. Actually, one can also see this directly
by using a pl-isomorphism of fans $\Delta $ and $\Phi \times \rho$, defined by
the function $f$.
\end{pf}

\section{Some immediate applications and generalizations}

\subsection{}
Let $\Phi $ be a complete fan in $V$, $\Delta \subset \Phi$ - a subfan which is quasi-convex.
 Then the fan $\Sigma :=\Phi -\Delta ^0$ is also quasi-convex, since
 $\partial \Delta =\partial \Sigma $. Consider the exact sequences of free $A$-modules
$$0\to \Gamma _{\Delta ^0}\cL _{\Phi}\to \Gamma (\Phi ,\cL _{\Phi})\to
\Gamma (\Sigma ,\cL _{\Phi})\to 0,$$
$$0\to \Gamma _{\Sigma ^0}\cL _{\Phi}\to \Gamma (\Phi ,\cL _{\Phi})\to
\Gamma (\Delta ,\cL _{\Phi})\to 0.$$
The pairings $[\cdot ,\cdot ]$ on the three fans identify the terms of the second sequence
as the duals (in the sense of $\Hom _A(\cdot ,A(2n)$) of the corresponding terms of the
first one.
It follows from Proposition 7.1 that the whole second sequence is the dual of the first one,
i.e. the maps in the second sequence are  adjoint to the maps in the first.

\subsection{} Let us translate Timorin's Theorem 4.7 to the fans setting. Given a
complete simplicial fan $\Phi$ choose a convex simple
$n$-dimensional polytope $P\subset W$
 so that $\Phi =\Phi _P$ . Let $P_i\subset P$ be one of the $n-1$-dimensional
faces and
$\rho _i\in \Phi$ -- the corresponding $1$-dimensional cone. Then the differential operator
$\partial _i\in A_1(P)$ corresponds to a piecewise linear function $\lambda _i$ on $\Phi$ that
takes nonzero values on $\rho _i$ and is zero on every other $1$-dimensional cone. Given a face
of $F=P_{i_1}\cap ...\cap P_{i_k}$ of $P$ the product $\lambda _{i_1}\cdot ...\cdot
\lambda _{i_k}$ is supported on the star $\St(\sigma _F)$.
As an immediate corollary of Theorem 4.7 and Proposition 7.1 we obtain the following proposition.

\begin{prop} For a complete simplicial fan $\Phi$ the space $H^{2k} (\Phi)=IH^{2k}(\Phi)$
has a basis consisting of (residues of) functions $f\in \cA (\Phi)$ such that the support
of $f$ is contained in $\St (\sigma )$ for $\sigma \in \Phi$ of dimension $k$. In particular,
the map
$$\bigoplus_{\dim (\sigma )=1}H^\bullet([\St (\sigma )],\partial \St (\sigma))\to
\bigoplus_{k>0}H^{2k}(\Phi)$$
is surjective. Hence the dual map
$$\bigoplus_{k<n}H^{2k}(\Phi)\to
\bigoplus_{\dim (\sigma )=1}H^\bullet([\St (\sigma )])$$
is injective.
\end{prop}

\subsection{} Following [Ka] we present a simple argument which deduces the (HL) theorem
for simplicial fans of dimension $n$ from the (HR) theorem for simplicial fans of
dimension $n-1$.

Let  $\Phi$ be a complete simplicial fan in $V$ and $l$ be a strictly convex
piecewise linear function on $\Phi$. After adding a globally linear function to it we may
assume that $l$ is strictly positive on $V-0$.
For each $1$-dimensional cone $\rho _i \in \Phi $ consider the $(n-1)$-dimensional complete
fan $\Phi _ i=p _ i(\partial \St(\rho _ i))$ in $V _ i =V/\langle \rho _i \rangle $
(where $p_ i:V\to V _ i$ is the projection).
Since the fan $\Phi$ is simplicial we   can write (uniquely)
$$l=\sum _ i \lambda _ i ,$$
where $\lambda _ i $ is a piecewise linear function on $\Phi$ supported on $\St (\rho _ i)$.

Suppose $h\in IH^{n-k}(\Phi)$ is such that $$l^k\cdot h =0.$$ Then
$$0=(h,l^kh)_{\Phi}=\sum _ i(h,\lambda _ i l^{k-1}h)_{\Phi}.$$

Fix a $1$-dimensional cone $\rho _i\in \Phi$.
 Let
$h _ i \in IH([\St (\rho _i)])$ denote the restriction of $h$ to
$\St (\rho _i )$. Changing $l$ by a global linear function (which depends on $i$) we may
assume that $l\vert _{\St (\rho _ i)}=p^* _i (l _ i)$ for a
strictly convex function $l _ i $ on $\Phi _ i$.  Identifying
$IH([\St (\rho _i)])=IH(\Phi _i)$ as in Lemma 3.18 we find that
$$h _ i \in Prim _{l_i}IH(\Phi _i).$$

Since $\lambda _i $ is supported on $\St (\rho _i )$ by
Proposition 7.1 $$(h,\lambda _ i l^{k-1}h)_{\Phi}=(h_i , \lambda
_i l^{ k-1} h _i)_{[\St (\rho _i )]}.$$ Now apply Proposition 7.8
to get
 $$(h_i , \lambda _i l^{ k-1} h _i)_{[\St (\rho _i)]}=
 n(h_i , l _i ^{ k-1} h _i)_{\Phi _i}.$$
By the (HR) theorem for $\Phi _i$, $(-1)^{\frac{n-k}{2}}(h_i , l^{
k-1} h _i)_{\Phi _i}>0$. It follows that $h _i=0$ for all $i$. But
then by Proposition 8.1 $h=0$.

This argument can be used to prove the following corollary

\begin{cor}[Ka] Let $\Phi$ be a complete $n$-dimensional
fan with a strictly convex function $l$.
Assume that $\Phi ^s\subset [\sigma ]$ for a unique cone $\sigma \in \Phi$.
Let $h\in IH^{n-k}(\Phi )$ be such that $l^kh=0$. Then the (HR) theorem in dimension $n-1$
implies that the restriction $h_i\in IH^{n-k}([\St (\rho _i)])$ is zero for all
$1$-dimensional cones $\rho _i\in \Phi -[\sigma]$.
\end{cor}

\begin{pf} Let $\rho _i$ be a $1$-dimensional cone in $\Phi -[\sigma ]$. Then as in the
previous argument (and using the same notation) we find that
$$h_i\in Prim _{l_i}IH^{n-k}(\Phi _i).$$ By adding a global linear
function to $l$ we may assume that $l\vert _{\sigma }=0$ and
$l(x)>0$ for $x \notin \sigma $. Then as before we can write
$$l=\sum _i\lambda _i,$$ where the summation is over all $\rho
_i\in \Phi -[\sigma ]$. Then again $$0=(h,l^kh)_{\Phi}=\sum _
i(h,\lambda _ i l^{k-1}h)_{\Phi}= n\sum _i
(h_i,l_i^{k-1}h_i)_{\Phi _i}.$$ If $h_i\neq 0$, then by (HR)
theorem for $\Phi _i$,
$(-1)^{\frac{n-k}{2}}(h_i,l_i^{k-1}h_i)_{\Phi _i}>0$. This proves
the corollary.
\end{pf}

\subsection{} Sometimes it is useful to have a pairing on $\Gamma
(\cL_\Delta)$ for a quasi-convex fan. The next proposition follows
easily from the results about the usual pairing. We denote by
$Q(A)$ the localization of $A$ with respect to all nonzero
homogeneous polynomials.

\begin{prop} Let $\Delta $ be a quasi-convex fan in $V$.  There exists a
unique $\cA (\Delta)$-bilinear pairing $$\{ \cdot ,\cdot \}
_{\Delta}:\Gamma (\cL _{\Delta})\times \Gamma (\cL _{\Delta})\to
Q(A)(2n),$$ which has the following properties.

a) The restriction of $\{ \cdot ,\cdot \} _{\Delta }$ to $\Gamma
(\cL _{\Delta})\times \Gamma _{\Delta ^0}\cL _{\Delta}$ takes
values in $A(2n)$ and coincides with $[\cdot ,\cdot ]_{\Delta}$.
In particular, if $\Delta $ is complete then $\{ \cdot, \cdot
\}_{\Delta }=[\cdot ,\cdot ]_{\Delta}$.

b) Assume that the fan $\Delta$ is covered by quasi-convex subfans
$\Delta _i$ so that every cone in $\Delta$ of maximal dimension
belongs to a unique $\Delta _i$. Then the restriction map $$\Gamma
(\cL _{\Delta})\to \bigoplus_i\Gamma (\cL _{\Delta _i})$$ is an
isometry. That is, if $a_i\in \Gamma (\cL _{\Delta _i})$ denotes
the image of $a \in \Gamma (\cL _{\Delta})$ then
$$\{a,b\}_{\Delta}=\sum _i\{a_i,b_i\}_{\Delta _i}.$$

c) Let $\pi :\Sigma \to \Delta$ be a subdivision. Then an admissible morphism
$\alpha :\cL _{\Delta}\to \pi_*\cL _{\Sigma}$ defines an isometry
$\alpha :\Gamma (\cL _{\Delta})\to \Gamma (\pi_*\cL _{\Sigma})$ with respect to
the pairings $\{ \cdot ,\cdot \}$.

d) Let $\Phi$ be a complete fan with a covering $\Phi =\cup
_i\Delta _i$ by quasi-convex subfans $\Delta _i$, such that every
cone in $\Phi$ of maximal dimension belongs to a unique $\Delta
_i$. Then the restriction map $$\Gamma (\cL _{\Phi})\to
\bigoplus_i\Gamma(\cL _{\Delta _i})$$ is an isometry. That is if
$a_i\in \Gamma (\cL _{\Delta _i})$ denotes the image of $a\in
\Gamma (\cL _{\Phi})$, then $$[a,b]_{\Phi}=\sum_i\{ a_i
,b_i\}_{\Delta _i}.$$
\end{prop}

\begin{pf} Let us first prove uniqueness. Indeed, since the pairing $[\cdot ,\cdot ]_{\Delta}$
is nondegenerate,  the free $A$-modules $\Gamma _{\Delta ^0}\cL
_{\Delta}$ and $\Gamma (\cL _{\Delta})$ have the same rank. Hence
the $A$-module $\Gamma (\cL _{\Delta})/\Gamma _{\Delta ^0}\cL
_{\Delta}$ is torsion. Since the $A$-module $Q(A)(2n)$ is torsion
free it follows that a pairing $\{ \cdot ,\cdot \}_{\Delta}$
satisfying a) is unique.

The construction of the pairing $\{ \cdot ,\cdot \}$ is the "same"
as that of $[\cdot ,\cdot]$. Assume first that the fan $\Delta $
is simplicial, i.e. $\cL _{\Delta} = \cA _{\Delta}$. Then put $$\{
f,g\}_{\Delta}=\frac{1}{n!}\zeta(fg)\in Q(A),$$ where $\zeta$ is
the Brion functional. As was explained above this definition is
forced on us if we want the property a) to hold. For a general
quasi-convex fan $\Delta$ choose a subdivision $\pi: \Sigma \to
\Delta$ with a simplicial $\Sigma$ and an admissible morphism
$\alpha :\cL _{\Delta}\to \pi_*\cL _{\Sigma}$. Put
$$\{a,b\}_{\Delta}:=\{\alpha (a),\alpha (b)\}_{\Sigma}.$$ Here
again we had no choice if the property c) is to hold. Theorem 7.2
implies that the property a) holds, which in turn means that the
pairing  $\{ \cdot ,\cdot \}_{\Delta}$ is independent of the
choice of the subdivision $\pi$ and the admissible morphism
$\alpha$.

The properties b) and c) are easy to check and d) is a special
case of b).
\end{pf}

\begin{example} In the above proposition one can take the "affine" fan $\Delta =[\sigma]$ for
a cone $\sigma$ of dimension $n$. This gives  the pairing on the
stalk at the closed point $$\{ \cdot ,\cdot \}_{[\sigma ]}:\cL
_{[\sigma ] ,\sigma}\times \cL _{[\sigma ] ,\sigma}\to Q(A)(2n).$$
\end{example}

\subsection{} Consider a complete fan $\Phi$ in $V$ with the following two
properties: 1) $\Phi$ contains a $1$-dimensional cone $\rho$ and
its negative $\rho ^\prime :=-\rho$, and it has a local product
structure at these cones; 2) The quasi-convex subfans $\Delta
^+=[\St (\rho)]$ and $\Delta ^-= [\St (\rho ^\prime)]$ intersect
along the common boundary and $\Phi =\Delta ^+\cup \Delta ^-$.

Let $p:V\to \overline{V}=V/\langle \rho \rangle$ be the projection
and $\Phi _{\rho}=p(\partial \St (\rho))$ - the complete fan in
$\overline{V}$. Put $B=\Sym \overline{V}^*$.
 By Lemma 3.18
the map $p$ induces natural isomorphisms
$$\Gamma (\cL _{\Delta ^-})\simeq A\otimes _B \Gamma (\cL _{\Phi _{\rho}})
\simeq \Gamma (\cL _{\Delta ^+}).$$
Hence we obtain the isomorphism of $A$-modules $\gamma :\Gamma (\cL _{\Delta ^-})
\to \Gamma (\cL _{\Delta ^+}).$

\begin{lemma} For $a,b\in \Gamma (\cL _{\Delta ^-})$ we have
$$\{ a,b\}_{\Delta ^-}=-\{\gamma (a),\gamma (b)\}_{\Delta ^+}.$$
\end{lemma}

\begin{pf} Let $\pi _{\rho}:\Psi _{\rho}\to \Phi _{\rho}$ be a subdivision with a
simplicial $\Psi _{\rho}$. This induces subdivisions $\pi
_-:\tilde{\Delta} ^-\to \Delta ^-$, $\pi _+:\tilde{\Delta} ^+\to
\Delta ^+$ with simplicial $\tilde{\Delta}^-$ and
$\tilde{\Delta}^+$. An admissible morphism $\alpha _{\rho}:\cL
_{\Phi _{\rho}}\to \pi _{\rho * } \cL _{\Psi _{\rho}}$ induces
corresponding admissible morphisms $\alpha _-$ and $\alpha _+$.
The induces morphisms $\alpha _-$ and $\alpha _+$ on global
sections commute with the isomorphism $\gamma$. Hence by
Proposition 8.3 c) we may assume that the fan $\Phi $ is
simplicial. But then the lemma follows from the definition of the
Brion functional $\zeta $.
\end{pf}

\section{Proof of Hodge-Riemann and Hard Lefschetz theorems}

We assume that (HL) and (HR) theorems hold for fans of dimension $\leq n-1$.

Let $\Phi$ be a complete projective fan of dimension $n$ in $V$.
Consider its singular subfan $\Phi ^s$. If $\Phi ^s$ is empty then
$\Phi$ is simplicial and the (HL) and (HR) theorems hold for
$\Phi$ (Corollary 6.4). Otherwise choose a maximal cone $\sigma
\in \Phi ^s$ and a ray $\rho \in \sigma ^0$ and consider the
corresponding star subdivision $$\pi :\Psi \to \Phi$$ (Example
2.5). The fan $\Psi$ is also projective (Lemma 2.13) and is "less"
singular than $\Phi$, i.e. $\Psi ^s$ contains a smaller number of
cones than $\Phi ^s$. So by induction on the size of the singular
subfan we may assume that (HL) and (HR) theorems hold for $\Psi$.
We are going to deduce from this that the theorems hold for
$\Phi$. Choose a strictly convex piecewise linear function $l$ on
$\Phi$.

Let $\tilde{l}$ be a piecewise function on $\Psi$ with support in
$\St _{\Psi}(\rho)$ such that $\tilde{l}\vert _{\rho}<0$. By Lemma
2.13 if we choose $\tilde{l}$ sufficiently small then
$\hat{l}=l+\tilde{l}$ is strictly convex on $\Psi$. Fix one such
$\hat{l}$.

Let $p:V\to V/\langle \rho \rangle =\overline{V}$ be the
projection. Then $p(\partial \St(\rho))=\Phi _{\rho}$ is a
complete fan in $\overline{V}$. By Lemma 3.18 the projection $p$
induces an isomorphism of $\cA (\Phi _{\rho})$-modules $IH(\Phi
_{\rho })=IH([\St (\rho)]).$  Put $B=\Sym \overline{V}^*$.

{\it Case 1:} Assume that $\dim (\sigma )=n$. We want to do this special case first,
because it follows almost immediately from the functorial properties of the pairing $(\cdot ,
\cdot )$, and the main ideas are transparent.

 Choose an admissible embedding $\alpha :\cL _{\Phi }\to \pi _*\cL _{\Psi}$. This induces
 an isometry $\alpha :IH(\Phi)\hookrightarrow IH(\Psi)$. We claim that one can choose
 $\alpha$ so that $\alpha (Prim _lIH(\Phi))\subset Prim _{\hat l }IH(\Psi)$. This is proved in
 Lemma 9.2 below.

Choose linear coordinates $x_1,...,x_n$ on $V$ so that the first
$n-1$ vanish on $\rho $ and $x_n$ is negative on $\St(\rho)$. Then
$\tilde{l} =x_n+f(x_1,...,x_{n-1})$ where $f$ is a strictly convex
(because the cone $\sigma $ is convex)
 piecewise linear function on $\Phi _{\rho }$. Therefore the $\tilde{l}$-action on
 $IH([\St (\rho)])$ coincides with the $f$-action on $IH(\Phi _{\rho })$, so,
 in particular,
 $$\tilde{l}^i:IH^{n-1-i}([\St (\rho)])\to IH^{n-1+i}([\St (\rho )])$$
 is an isomorphism.

\begin{lemma} We may choose an admissible embedding $\alpha :\cL _{\Phi}\hookrightarrow
\pi_*\cL _{\Psi}$ so that the image of $\overline{\cL _{\Phi ,\sigma}}$ in
$IH([\St (\rho)])$ is equal to the
$\tilde{l}$-primitive subspace.
\end{lemma}

\begin{pf} This is clear from the definition of the sheaf $\cL
_{\Phi}$.
\end{pf}

\begin{lemma} Under an admissible embedding $\alpha $ as in the last lemma the
image of $Prim_l IH(\Phi)$ is contained in $Prim _{\tilde{l}}IH(\Psi)$.
\end{lemma}

\begin{pf} For simplicity of notation we identify $IH(\Phi)$ as a
subspace of $IH(\Psi)$ (by means of $\alpha $).

Choose $a\in Prim_lIH^{n-k}(\Phi)$. We have
$$\hat{l}^{k+1}a=l^{k+1}a+\tilde{l}q(l,\tilde{l})a
=\tilde{l}q(l,\tilde{l})a$$ for a polynomial $q$. Since the
support of $\tilde{l}$ is contained in $\St(\rho)$ and
$l\vert_{\St(\rho)}$ is linear we have $\tilde{l}la=0$. Thus
$$\hat{l}^{k+1}a=\tilde{l}^{k+1}a$$ may be considered as an
element of $IH([\St (\rho)],\partial \St (\rho))$ (which is a
subspace of $IH(\Psi)$). We have   $\tilde{l}^{k+1}a=0$ if
$(c,\tilde{l}^{k+1}a)_{[\St (\rho )]}=0$ for all $c\in IH([\St
(\rho)]$. By Proposition 7.8 $$(c,\tilde{l}^{k+1}a)_{[\St (\rho)]}
=
 -n(c,f^ka)_{\Phi _{\rho}}=0,$$ because
 $f^ka=0$.
\end{pf}

Let us fix an admissible embedding $\alpha $ as in Lemma 9.1 above
and choose $0\neq a\in Prim _lIH^{n-k}(\Phi)$. Then by Lemma 9.2
and by our induction hypothesis
$$(-1)^{\frac{n-k}{2}}(a,\hat{l}^ka)_{\Psi}>0.$$ We have
\begin{eqnarray}
(a,\hat{l}^ka)_{\Psi} & = & (a,l^ka)_{\Psi}+
(a,\tilde{l}t(l,\tilde{l})a)_{\Psi}\nonumber\\ & = &
(a,l^ka)_{\Phi}+ (a,\tilde{l}t(l,\tilde{l})a)_{\Psi}\nonumber
\end{eqnarray}
for a polynomial $t$. By the support consideration the element
$\tilde{l}t(l,\tilde{l})a$ belongs to the subspace $IH([\St
(\sigma )],\partial \St (\sigma ))\subset IH(\Psi)$. By abuse of
notation we will denote also by $a$ the image of $a$ in $IH([\St
(\sigma )])$. Then by Propositions 7.1 and 7.8
\begin{eqnarray}
(a,\tilde{l}t(l,\tilde{l})a)_{\Psi} & = &
(a,\tilde{l}t(l,\tilde{l})a)_{[\St (\sigma )]}\nonumber\\ & = &
(a,\tilde{l}^ka)_{[\St (\sigma )]}\nonumber\\ & = &
-n(a,f^{k-1}a)_{\Phi _{\rho}} \nonumber
\end{eqnarray}

By Lemma 9.1 the element $a\in IH^{n-k}(\Phi _{\rho})$ is
$f$-primitive. Hence by (HR) theorem for $\Phi _{\rho}$
$$(-1)^{\frac{n-k}{2}} (a,f^{k-1}a)_{\Phi _{\rho}}\geq 0$$ (the
inequality is not strict since $a$ may be zero in  $IH^{n-k}(\Phi
_{\rho})$). Therefore $$ (a,\tilde{l}t(l,\tilde{l})a)_{\Psi}\leq
0$$ and $$ (-1)^{\frac{n-k}{2}}(a,l^ka)_{\Phi}>0,$$ which proves
the (HL) and (HR) theorems for $\Phi$.

\begin{remark} Our proof shows that the quadratic form
$Q_{\hat{l}}$ on $IH(\Psi)$ tends to be more degenerate than the
form $Q_l$ on $IH(\Phi)$. Indeed, this is what happens in blow-ups
 of algebraic varieties. We will also see this in the remaining case
 below.
\end{remark}

{\it Case 2:} $\dim (\sigma )<n$. Here we essentially copy [Ka].

Denote $\Delta =[\St _{\Phi}(\sigma )]$, $\rho ^\prime =-\rho$. Consider
the "complementary" fan $$\Delta ^\prime :=\{ \rho ^\prime +\tau
\vert \tau \in \partial \Delta \}.$$ The (quasi-convex) fans
$\Delta $ and $\Delta ^\prime $ intersect along the common
boundary and their union is a complete fan in $V$ which we denote
$\overline{\Delta }$. First we prove the (HR) theorem for this
auxiliary fan $\overline{\Delta}$ and then use a trick to deduce
the (HR) theorem for $\Phi$.

Denote $\hat{\Delta}=[\St _{\Psi}(\rho)]$
(so that $\vert \Delta \vert =\vert \hat{\Delta } \vert$).

By changing $\hat{l}$ by a global linear function we may assume
that $\hat{l}\vert _{\rho}=0$ and hence $\hat{l}\vert
_{\hat{\Delta}}=p^*l_{\rho}$ for a strictly convex function
$l_{\rho}$ on $\Phi _{\rho}$. Define a piecewise linear function
$\bar{l}$ on the fan $\overline{\Delta}$ as follows
$$\bar{l}(x)=\begin{cases}l(x) & \text{if $x\in \vert \Delta
\vert$}\\ l_{\rho}(p(x)) & \text{if $x\in \vert \Delta ^\prime
\vert .$}
\end{cases}
$$ Then $\bar{l}$ is strictly convex on $\overline{\Delta}$.
Indeed, this follows from the inequality
$l(x)>\hat{l}(x)=l_{\rho}(p(x))$ for $x$ in the interior of $\vert
\Delta \vert$.

\begin{prop} The map
$$\bar{l}^k:IH^{n-k}(\overline{\Delta})\to
IH^{n+k}(\overline{\Delta})$$ is an isomorphism for all $k\geq 1$.
That is the (HL) theorem holds for the fan $\overline{\Delta}$
with the strictly convex function $\bar{l}$.
\end{prop}

\begin{pf} Let $h\in IH^{n-k}(\overline{\Delta})$ be in the kernel
of $\bar{l}^k$.

Consider the exact sequence of free $A$-modules $$0\to \Gamma
_{\Delta ^0}\cL _{\overline{\Delta}}\to \Gamma (\cL
_{\overline{\Delta}})\to \Gamma (\Delta ^\prime ,\cL
_{\overline{\Delta}})\to 0$$ and the induced exact sequence of
graded vector spaces $$0\to IH(\Delta ,\partial \Delta
 )\to IH(\overline{\Delta})\to IH(\Delta ^\prime)\to 0.$$ Note that
 $\overline{\Delta}^s\subset [\sigma]\subset \Delta$. Hence by
 Corollary 8.2, $h\in IH(\Delta ,\partial \Delta)$. Since multiplication
 by $\bar{l}$ is a self adjoint operator with respect to the
 pairing $(\cdot ,\cdot)_{\overline{\Delta}}$ it suffices to prove
 the following lemma.

\begin{lemma} The map
$$l^k:IH^{n-k}(\Delta)\to IH^{n+k}(\Delta)$$ is surjective for all
$k\geq 1$.
\end{lemma}

\noindent{\it Proof of lemma.} Let $d=\dim(\sigma)$ and consider
the projection $r:V\to V/\langle \sigma \rangle$. Then $\Phi
_{\sigma}=r(\Link (\sigma))$ is a complete fan in $V/\langle
\sigma \rangle$. Changing $l$ by a global linear function we may
assume that $l\vert _\Delta=r^*(l_{\sigma })$ for a strictly
convex function $l_{\sigma }$ on $\Phi _{\sigma }$. Then Lemma
3.18 implies that $$IH(\Delta)\simeq IH(\Phi _{\sigma })\otimes
\overline{\cL _{\Delta ,\sigma }},$$ where $l$ acts as
$l_{\sigma}\otimes id$. The assertion of the lemma follows
immediately from the (HL) theorem for $\Phi _{\sigma }$ and the
fact that the graded vector space $\overline{\cL _{\Delta ,\sigma
}}$ is zero in degrees $\geq d$ (3.13). This proves the lemma and
the proposition.
\end{pf}

\begin{lemma} The quadratic form $Q_{\bar{l}}$ on $IH(\overline{\Delta})$ satisfies the
Hodge-Riemann bilinear relations.
\end{lemma}

\begin{pf}
Choose a subspace $V_1\subset V$ complementary to $\langle \sigma
\rangle$ and identify $V_1=V/\langle \sigma \rangle$ by the
projection $r$, so that $\Phi _{\sigma }$ is a complete fan in
$V_1$. Denote by $\overline{[\sigma]}\subset \overline{\Delta}$
the subfan lying in the subspace $\langle \sigma \rangle$. Then
$\overline{[\sigma]}$ is a complete subfan in $\langle \sigma
\rangle$, whose $1$-dimensional cones are those of $\sigma $ and
$\rho ^\prime$.

The projection $r:V\to V_1$ defines a pl-isomorphism of fans
$\overline{\Delta}\to \Phi _{\sigma}\times \overline{[\sigma]}$.
Let $\phi : \Phi _{\sigma}\times \overline{[\sigma]}\to
\overline{\Delta}$ be the inverse isomorphism. Then $\phi $ is
determined by a function $g:V_1\to \langle \sigma \rangle$ which
is piecewise linear with respect to the fan $\Phi _{\sigma}$ ($\Link
(\sigma )$ is the graph of $g$). For each $0\leq t\leq 1$ the
function $tg$ defines similarly a pl-isomorphism $$\phi _t: \Phi
_{\sigma}\times \overline{[\sigma]}\to \overline{\Delta}_t,$$
Where $\overline{\Delta}_1=\overline{\Delta}$ and
$\overline{\Delta}_0= \Phi _{\sigma}\times \overline{[\sigma]}$.
Thus $\Phi _{\sigma}\times \overline{[\sigma]}$ and
$\overline{\Delta}$ are two members of a continuous family of fans
$\overline{\Delta}_t$.

As in the proof of last lemma we can assume that $\bar{l}\vert
_{\sigma }=0$ and hence $\bar{l}\vert _\Delta =r^*l_{\sigma }$ for
a strictly convex function $l_{\sigma }$ on $\Phi _{\sigma}$. We
write $$\bar{l}=r^*l_{\sigma }+(\bar{l}-r^*l_{\sigma }),$$ where
$(\bar{l}-r^*l_{\sigma })$ is supported on $\Delta ^\prime$. Note
that the restriction of $(\bar{l}-r^*l_{\sigma })$ to
$\overline{[\sigma ]}$ is strictly convex. Consider the piecewise
linear function $$\phi ^*\bar{l}=\phi ^*r^*l_{\sigma}+\phi
^*(\bar{l}-r^*l_{\sigma })$$ on $\Phi _{\sigma}\times
\overline{[\sigma]}$. The two summands are the pullbacks of
strictly convex functions from $\Phi _{\sigma }$ and
$\overline{[\sigma ]}$ respectively. For each $0\leq t\leq 1$
consider the piecewise linear function $$\bar{l}_t:=\phi _t^{-1
*}\phi ^*\bar{l}.$$ This function is strictly convex on
$\overline{\Delta}_t$. Repeating the proof of Proposition 9.4 we
find that (HL) theorem holds for the operators $\bar{l}_t$ on
$IH(\overline{\Delta }_t)$. Thus we obtain a continuous family of
graded vector spaces $IH(\overline{\Delta }_t)$ with
continuously varying operators $\bar{l}_t$. By the (HL)-property
the corresponding quadratic forms $Q_{\bar{l}_t}$ have full rank.
By Theorem 10.8  the form $Q_{\bar{l}_0}$ satisfies the
Hodge-Riemann relations. Hence so do all the forms
$Q_{\bar{l}_t}$. This proves the (HR) theorem for
$\overline{\Delta}$.
\end{pf}

To complete the proof of (HR) theorem for $\Phi$ we will deduce it from the
same theorem for $\Psi $ and $\overline{\Delta}$.

The fans $\Phi $ and $\overline{\Delta }$ have a common subfan
$\Delta$. Define the $A$-module $$F:=\{ (s_1,s_2)\in \Gamma (\cL
_{\Phi})\times \Gamma ( \cL _{\overline{\Delta}})\vert \ \ \
s_1\vert _\Delta =s_2\vert _\Delta \}$$ with the operator
$l_F=(l,\bar{l})$ and the bilinear form $$[\cdot ,\cdot
]_F:F\times F\to A(2n),\quad [\cdot ,\cdot ]_F=[\cdot ,\cdot
]_\Phi - [\cdot ,\cdot ]_{\overline{\Delta}}.$$

\begin{lemma} There exists a morphism of $A$-modules $\beta :F\to \Gamma (\cL _{\Psi})$,
such that

a) $\beta \cdot l_F=\hat{l}\cdot \beta$;

b) $[a,b]_F=[\beta (a),\beta (b)]_{\Psi}$.
\end{lemma}

\begin{pf}
By Lemma 3.18 applied to the fans $\hat{\Delta}$ and $\Delta ^\prime$ and the projection
$p:V\to \overline{V}$ we obtain natural isomorphisms of $A$-modules (also of
$\Gamma (\cL _{\Phi _{\rho}})$-modules)
$$\Gamma (\cL _{\Delta ^\prime})\simeq A\otimes _B \Gamma (\cL _{\Phi _{\rho}})
\simeq \Gamma (\cL _{\hat{\Delta}}).$$
This defines a natural isomorphism $\gamma :\Gamma (\cL _{\Delta ^\prime})\to
\Gamma (\cL _{\hat{\Delta}})$.

Consider the exact sequence of $A$-modules
$$0\to \Gamma (\cL _{\Psi})\to
\Gamma (\Psi -\hat{\Delta }^0 ,
\cL _{\Psi})\oplus \Gamma (\hat{\Delta} ,\cL _{\Psi})\stackrel{(+,-)}
{\rightarrow} \Gamma (\partial \hat{\Delta} ,\cL _{\Psi})\to 0.$$

Define $\beta :\Gamma (\cL _{\Phi})\oplus \Gamma (\cL _{\overline{\Delta}})\to
\Gamma (\Psi -\hat{\Delta }^0 ,
\cL _{\Psi})\oplus \Gamma (\cL _{\hat{\Delta}})$ to be
$$\beta (s_1,s_2):=(s_1\vert _{\Phi -\Delta ^0}, \gamma (s_2\vert _{\Delta ^\prime})).$$
This map descends to a map
$$\beta :F\to \Gamma (\cL _{\Psi}).$$
This is a map of $A$-modules which satisfies  the property a) by the choice of functions
$\bar{l}$ and $\hat{l}$. Let us show that it satisfies b).

Consider the composition of natural maps
$$\begin{array}{rcl}
 & F & \\
 & \downarrow & \\
\Gamma (\cL _{\Phi}) & \oplus &  \Gamma (\cL _{\overline{\Delta}})\\
\downarrow & & \downarrow\\
\Gamma (\Phi -\Delta ^0,\cL _{\Phi}) \oplus \Gamma (\Delta ,\cL _{\Phi}) & \oplus &
\Gamma (\Delta ^\prime ,\cL _{\overline{\Delta}})\oplus
\Gamma (\Delta ,\cL _{\overline{\Delta}})\\
(id,0)\downarrow & & \downarrow (\gamma ,0)\\
\Gamma (\Psi -\Delta ^0,\cL _{\Psi}) & \oplus & \Gamma (\hat{\Delta},\cL _{\Psi}).
\end{array}$$
The image of $F$ under this composition is contained in $\Gamma (\cL_{\Psi})$ (and
the composition itself is equal to $\beta$). It is convenient to use the pairing
$\{ \cdot ,\cdot \}$ (8.3) on each summand in this diagram. More precisely, on the
summands involving $\cL _{\overline{\Delta}}$ we take the pairing $\{ \cdot ,\cdot \}$
with the negative sign. It then suffices to prove that the composition of maps is
an isometry. The first map is such by definition. The second is by Proposition 8.3.
The third map is not an isometry, but it is such on the image of $F$. This proves the
lemma.
\end{pf}

Let $\overline{F}$ as usual denote the graded vector space $F/A^+F$.
Denote by
$$(\cdot ,\cdot )_F:\overline{F}\times \overline{F}\to \bbR (2n)$$
the
residue of the pairing $[\cdot ,\cdot ]_F$.
There is a natural
 map $\overline{F}\to IH(\Phi)\oplus IH(\overline{\Delta})$.

\begin{lemma} Let $a\in Prim _lIH^{n-k}(\Phi)$. Then there exists $b\in Prim_{\bar{l}}
IH^{n-k}(\overline{\Delta})$ such that $(a,b)$ belongs to the image of $\overline{F}$.
\end{lemma}

\begin{pf} Denote by $c\in IH^{n-k}(\Delta)$ the image of $a$.
Then $l^{k+1}c=0$. We need to find $b\in
Prim _{\bar{l}}IH^{n-k}(\overline{\Delta})$ such that $b\vert _{\Delta}=c$.
Consider the exact sequence
$$0\to IH(\Delta ^\prime ,\partial \Delta ^\prime)\to IH(\overline{\Delta})\to
IH(\Delta)\to 0.$$
Let $b^\prime \in IH^{n-k}(\overline{\Delta})$ be any preimage of $c$. Then
$\bar{l}^{k+1}b^\prime \in IH^{n+k+2}(\Delta ^\prime ,\partial \Delta ^\prime)$.
Note that the action of the operator $\bar{l}\vert _{\Delta ^\prime}$ on $IH(\Delta ^\prime)$
coincides with the action of $l_{\rho}$ on $IH(\Phi _{\rho})$ (3.18). Hence by the induction
hypothesis for $\Phi _{\rho}$ the map
$$\bar{l}^{k+1}:IH^{n-k-2}(\Delta ^\prime)\to IH^{n+k}(\Delta ^\prime)$$
is an isomorphism. By Corollary 7.9 the $\cA (\Delta ^\prime)$-modules
$IH(\Delta ^\prime)$ and $IH(\Delta ^\prime ,\partial \Delta ^\prime)$ are isomorphic with a
shift by 2. Hence the map
$$\bar{l}^{k+1}:IH^{n-k}(\Delta ^\prime ,\partial \Delta ^\prime)\to
IH^{n+k+2}(\Delta ^\prime ,\partial \Delta ^\prime)$$
is an isomorphism. So there exists $h \in
IH^{n-k}(\Delta ^\prime ,\partial \Delta ^\prime)$ such that $\bar{l}^{k+1}h=
\bar{l}^{k+1}b^\prime $. Now put $b=b^\prime -h$.
\end{pf}

Now the (HR) theorem for the fan $\Phi$ follows. Indeed, take $a\in Prim _l IH^{n-k}(\Phi)$.
 Let $b\in Prim_{\bar{l}}IH^{n-k}(\overline{\Delta})$ be as in last lemma. Then by Lemma
 9.7
 $$(-1)^{\frac{n-k}{2}}(a,l^ka)_{\Phi}-(-1)^{\frac{n-k}{2}}(b,\bar{l}^kb)_{\overline{\Delta}}$$
 $$=(-1)^{\frac{n-k}{2}}((a,b),l_F(a,b))_F=(-1)^{\frac{n-k}{2}}(\beta (a,b),\hat{l}\beta (a,b))_{\Psi}.$$
By the (HR) theorem for $\overline{\Delta}$ and $\Psi$ it follows
that $$(-1)^{\frac{n-k}{2}}(a,l^ka)_{\Phi}\geq 0.$$ Assume that
$a\neq 0$.
 If  $b\neq 0$ then
$$(-1)^{\frac{n-k}{2}}(a,l^ka)_{\Phi}\geq
(-1)^{\frac{n-k}{2}}(b,\bar{l}^kb)_{\overline{\Delta}}>0$$ and we
are done. Assume that $b=0$. Then $a\vert _{\Delta}=0$, and hence
$a\in IH(\Phi -\Delta ^0)$, so that $\beta (a,b)\neq 0$. Thus
$$(-1)^{\frac{n-k}{2}}(a,l^ka)_{\Phi}= (-1)^{\frac{n-k}{2}}(\beta
(a,b),\hat{l}\beta (a,b))_{\Psi}>0.$$ This completes the proof.

\section{K\"unneth formula for $IH$ and duality on the product of fans}

\subsection{K\"unneth formula}

Let $V_1$, $V_2$ be vector spaces of dimensions $n_1$ and $n_2$ respectively.
Let $A_1=\Sym V^*_1$, $A_2=\Sym V^*_2$ be the (evenly graded) algebras of polynomial functions
on $V_1$ and $V_2$ respectively. Put $V=V_1\times V_2$, $A=\Sym V^*$; then
$A=A_1\otimes A_2$. Denote by $V_1\stackrel{p_1}
{\leftarrow }V\stackrel{p_2}{\rightarrow}V_2$ the two projections.

Let $\Delta$ and $\Sigma $ be fans in $V_1$ and $V_2$ respectively. Consider the
{\it product}
fan $\Phi =\Delta \times \Sigma $ in $V$
$$\Phi =\{ \sigma +\tau \vert \sigma \in \Delta ,\tau \in \Sigma \}$$
with the projections $p_1:\Phi \to \Delta$, $p_2:\Phi \to \Sigma$.

Let  $F_1\in \Sh (\Delta)$, $F_2\in \Sh(\Sigma )$. Then for every $\sigma \in \Delta$
$$(p_{1*}(p^{-1}_1F_1\otimes p_2^{-1}F_2))_{\sigma }=F_{1,\sigma }\otimes \Gamma (\Sigma ,F_2).$$
In particular,
$$\Gamma (\Phi ,p_1^{-1}F_1\otimes p_2^{-1}F_2)=\Gamma (\Delta ,F_1\otimes \Gamma (\Sigma ,
F_2))=\Gamma (\Delta ,F_1)\otimes \Gamma (\Sigma ,F_2).$$

Note the  canonical
isomorphisms of sheaves
$$A_{\Phi}=p^{-1}_1A_{\Delta}\otimes p^{-1}_2A_{\Sigma }, \quad
\cA_{\Phi}=p^{-1}_1\cA_{\Delta}\otimes p^{-1}_2\cA_{\Sigma }.$$
In particular we have $\cA (\Phi )=\cA (\Delta )\otimes \cA (\Sigma )$.

Consider the $\cA _{\Phi}$-module $\cL _{\Phi}^\prime:=p_1^{-1}\cL _{\Delta}\otimes
p_2^{-1}\cL _{\Sigma }$.

\begin{lemma} There exists a canonical isomorphism
of $\cA _{\Phi}$-modules $\cL _{\Phi}^\prime = \cL _{\Phi}$, hence an isomorphism
$\Gamma (\Phi ,\cL _{\Phi})=\Gamma (\Delta ,\cL _{\Delta })\otimes \Gamma (\Sigma ,
\cL _{\Sigma }).$ Thus if $\Delta $ and $\Sigma $ are quasi-convex, then so is
$\Phi$, and $IH(\Phi)=IH(\Delta)\otimes IH(\Sigma )$.
\end{lemma}

\begin{pf} Clearly, $\cL _{\Phi ,\Or}^\prime =\bbR$ and the
 $\cA _{\Phi}$-module $\cL _{\Phi}^\prime$ is locally free. It suffices to show that
 the map $\overline{\cL _{\Phi ,\xi}^\prime}\to
 \overline {\Gamma (\partial \xi ,\cL^\prime _{\Phi})}$
 is an isomorphism for $\xi \in \Phi$.

 Let  $\xi=\sigma +\tau \in \Phi$ where $\sigma \in \Delta$, $\tau \in \Sigma$. Then
 $\partial \xi=[\sigma ]\times \partial \tau \cup \partial \sigma \times [\tau ]$
 and
 $$\Gamma (\partial \xi ,\cL ^\prime_{\Phi})=\Ker \{\Gamma
 ([\sigma ]\times \partial \tau ,\cL _{\Phi}^\prime )\oplus
 \Gamma (\partial \sigma \times [\tau ],\cL _{\Phi }^\prime)
 \stackrel{(+,-)}{\rightarrow}\Gamma (\partial \sigma \times \partial \tau ,
 \cL ^\prime _{\Phi})\}.$$
 Note that the complex
$$\Gamma ([\xi] ,\cL ^\prime_{\Phi})\rightarrow \Gamma
 ([\sigma ]\times \partial \tau ,\cL _{\Phi}^\prime )\oplus
 \Gamma (\partial \sigma \times [\tau ],\cL _{\Phi }^\prime)
 \stackrel{(+,-)}{\rightarrow}\Gamma (\partial \sigma \times \partial \tau ,
 \cL ^\prime _{\Phi})$$
is isomorphic to the tensor product of complexes
$$(*)\ \ \ \Gamma ([\sigma ],\cL _{\Delta})\to \Gamma (\partial \sigma ,\cL _{\Delta}),\quad
\Gamma ([\tau ],\cL _{\Sigma})\to \Gamma (\partial \tau ,\cL _{\Sigma }).$$
It follows that the map
$$\Gamma ([\xi] ,\cL ^\prime_{\Phi})\rightarrow \Gamma (\partial \xi ,\cL ^\prime_{\Phi})$$
is surjective (i.e. the sheaf $\cL ^\prime _{\Phi}$ is flabby), hence also the map
$\overline{\cL _{\Phi ,\xi}^\prime}\to
 \overline {\Gamma (\partial \xi ,\cL^\prime _{\Phi})}$ is such.
The two complexes $(*)$ (and so their tensor product)
become acyclic after taking the residue at the maximal ideal $A^+\subset A$.
It follows that the map
 $\overline{\cL _{\Phi ,\xi}^\prime}\to \overline {\Gamma (\partial ,\cL ^\prime_{\Phi})}$
 is also injective.
 \end{pf}

\begin{lemma} Assume that the fans $\Delta $ and $\Sigma $ (hence also $\Phi$) are quasi-convex.
Then
$$\Gamma _{\Phi ^0}\cL _{\Phi}=\Gamma _{\Delta ^0}\cL _{\Delta }\otimes
\Gamma _{\Sigma ^0} \cL _{\Sigma }.$$
Hence $IH(\Phi ,\partial \Phi)=IH(\Delta ,\partial \Delta)\otimes IH(\Sigma ,\partial
\Sigma )$.
\end{lemma}

\begin{pf}
By taking the tensor product of short exact sequences
$$0\to \Gamma _{\Delta ^0}\cL _{\Delta}\to \Gamma (\cL _{\Delta})\to
\Gamma (\partial \Delta, \cL _{\Delta})\to 0,$$
and
$$0\to \Gamma _{\Sigma ^0}\cL _{\Sigma}\to \Gamma (\cL _{\Sigma})\to
\Gamma (\partial \Sigma, \cL _{\Sigma})\to 0$$
we find the exact sequence
$$0\to \Gamma _{\Delta ^0}\cL _{\Delta}\otimes
\Gamma _{\Sigma ^0}\cL _{\Sigma }
\to \Gamma (\cL _{\Delta})
\otimes \Gamma (\cL _{\Sigma})$$
$$\to
[\Gamma (\partial \Delta, \cL _{\Delta})\otimes
\Gamma (\cL _{\Sigma })]\oplus
[\Gamma (\cL _{\Delta})\otimes
\Gamma (\partial \Sigma ,\cL _{\Sigma })].$$
Note that the kernel of the map
$$\Gamma (\cL _{\Delta})
\otimes \Gamma (\cL _{\Sigma})
\to [\Gamma (\partial \Delta, \cL _{\Delta})\otimes
\Gamma (\cL _{\Sigma })]\oplus
[\Gamma (\cL _{\Delta})\otimes
\Gamma (\partial \Sigma ,\cL _{\Sigma })]$$
is naturally isomorphic to the kernel of the map
$$\Gamma (\cL _{\Phi})\to \Gamma (\partial \Phi ,\cL _{\Phi}).$$
Hence
$$\Gamma _{\Phi ^0}\cL _{\Phi}=\Gamma _{\Delta ^0}\cL _{\Delta }\otimes
\Gamma _{\Sigma ^0} \cL _{\Sigma }.$$
\end{pf}

\subsection{Duality on the product of fans}

Consider the dualizing modules $\omega _1=A_1\cdot \wedge ^{n_1}V_1^*$,
$\omega _2=A_2\cdot \wedge ^{n_2}V_2^*$ on $V_1$ and $V_2$. Then $\omega =\omega _1\otimes
\omega _2$.

Choose volume forms $\Omega _{V_i}\in \wedge ^{n_i}V^*_{i}$, $i=1,2$ and put
$\Omega _V=\Omega _{V_1}\wedge \Omega _{V_2}$.

Assume that orientations of cones in $\Delta $ and $\Sigma $ are chosen (the cones
of top dimension are oriented by $\Omega _{V_1}$ and $\Omega _{V_2}$ respectively).
This determines the cellular complexes $C^\bullet _{\Delta}(\cdot)$ and
$C^\bullet _{\Sigma}(\cdot)$.
Given $\xi =\sigma + \tau \in \Phi$ with $\sigma \in \Delta$,
$\tau \in \Sigma$ we orient it in the usual way by putting first the
vectors in $\sigma $ and then the vectors in $\tau$. Thus the cones of dimension $n$ in
$\Phi$ are oriented by $\Omega _V$.
This defines the cellular complex $C^\bullet _{\Phi}(\cdot)$. Note that
$C^\bullet _{\Phi}=C^\bullet_{\Delta}\otimes C^\bullet _{\Sigma }$ as graded vector
spaces, but not as complexes.
I.e. for $F_1\in \Sh (\Delta)$, $F_2\in \Sh(\Sigma )$, $F=p_1^{-1}F_1\otimes
p_2^{-1}F_2\in \Sh(\Phi)$ we have
$$C^k_{\Phi}(F)=\bigoplus_s C^s_{\Delta}(F_1)\otimes C^{k-s}_{\Sigma }(F_2)$$
and the differential on the right is
$$d(a\otimes b)=d(a)\otimes b +(-1)^sa\otimes d(b),\quad \text{for $a\otimes b\in
C^s_\Delta(F_1)\otimes C^t_\Sigma (F_2)$},$$
whereas on the left it is
$$d(a\otimes b)=d(a)\otimes b +(-1)^{n_1+s}a\otimes d(b).$$
To fix this discrepancy let us consider the new cellular complex
$\overline{C}^\bullet$ which is equal to the shifted (to the left) complex $C^\bullet[d]$,
except it has the {\it same} differential (and not $(-1)^d$ times the original one), where
$d$ is the dimension of the ambient space. Then the identification
$$\overline{C}^\bullet _{\Phi}=\overline{C}^\bullet_{\Delta}\otimes
\overline{C}^\bullet _{\Sigma }$$
is the equality of complexes.

Note that we can use the complex $\overline{C}^\bullet$ instead of $C^\bullet$ in
Proposition 3.9. Namely for $F\in D^b_c(A_{\Phi}-mod)$ the functor $\overline{C}^\bullet
_{\Phi}$ induces the isomorphism
$$\bbR \Gamma (\Phi ,DF)\simeq \bbR \Hom ^\bullet _A(\overline{C}^\bullet,\omega [n]),$$
and similarly for $\Delta $ and $\Sigma $. We will use this description of the duality
to show that it is compatible with the product of fans.

Let $P_1^\bullet \in D^b_c(A_1-mod)$,  $P_2^\bullet \in D^b_c(A_2-mod)$,
 $P^\bullet =P_1^\bullet \otimes P_2^\bullet \in D^b_c(A-mod)$. Then there exists
 a natural functorial isomorphism in $D_c^b(A-mod)$:
$$\delta:\bbR \Hom _{A_1}(P_1^\bullet ,\omega _1[n_1])\otimes
\bbR \Hom _{A_2}(P_2^\bullet ,\omega _2[n_2])\to \bbR \Hom _A(P^\bullet ,\omega [n] ),$$
$$\delta(f\otimes g)(a\otimes b)=(-1)^{\deg (g)\deg (a)}f(a)\otimes g(b).$$

\begin{example} In the previous notation let $P_1^\bullet =\bbR = P_2^\bullet.$ As was
remarked in the proof of Lemma 3.15 above we have canonical isomorphisms
$\Ext _A^n(\bbR ,\omega [n])=\bbR =
\Ext _{A_i}^{n_i}(\bbR  ,\omega _i[n_i])$ for $i=1,2$. The map $\delta $ induces an
isomorphism
$$\delta : \Ext _{A_1}^{n_1}(\bbR  ,\omega _1 [n_1])\otimes
\Ext _{A_2}^{n_2}(\bbR ,\omega _2 [n_2])\to \Ext _A^n(\bbR  ,\omega [n])$$
which coincides with the multiplication map $\bbR \otimes \bbR \to \bbR$ under the
above isomorphisms. (Indeed, the tensor product of the Koszul resolutions of $\bbR$
as $A_1$- and $A_2$-module respectively is equal to its Koszul resolution as an $A$-module.)
\end{example}

\begin{prop} a) Let $F_1\in D^b_c(A_{\Delta}-mod)$, $F_2\in D^b_c(A_{\Sigma}-mod)$ and
$F=p_1^{-1}F_1\otimes p_2^{-1}F_2 \in D_C^b(A_{\Phi}-mod)$.
The map $\delta $ defines a functorial isomorphism in $D_c^b(A_{\Phi}-mod)$
$$\delta :  p_1^{-1}DF_1\otimes p_2^{-1}DF_2 \to DF.$$

b) In case $F_1=\cL _{\Delta}$, $F_2=\cL _{\Sigma}$ and $F=\cL _{\Phi}$
this isomorphism has the following properties:

i) If we use
the canonical identifications of stalks
$D(\cL)_{\Or}=\bbR$ on each of the three fans,
then $\delta _{\Or}$ is the multiplication map $\bbR \otimes \bbR \to \bbR$.

ii)  Assume that the fans $\Delta$, $\Sigma $ (and hence $\Phi$) are quasi-convex. We use
isomorphisms of Proposition 3.9 (with $\overline{C}^\bullet$ instead of $C^\bullet$)
and Lemma 10.2.
Then on the level of global sections the map
$$\delta :\Hom _{A_1}(\Gamma _{\Delta ^0}\cL _{\Delta},\omega _1)\otimes
\Hom _{A_2}(\Gamma _{\Sigma ^0}\cL _{\Sigma},\omega _2)\to
\Hom _A(\Gamma _{\Phi ^0}\cL _{\Phi},\omega )$$
is given by
$$\delta (f\otimes g)(a\otimes b)=f(a)\otimes g(b).$$
\end{prop}

\begin{pf} a) Let $\xi =\sigma +\tau \in \Phi$. It follows from Proposition 3.9 that
$DF([\xi])=
\bbR \Hom _A(\overline{C}_{\Phi}^\bullet (F_{[\xi]}),\omega [n])$ and
similarly for $DF_1([\sigma ])$ and
$DF_2([\tau ])$. Applying $\delta $ to the complexes $\overline{C}^\bullet _{\Phi }(F_{[\xi]})=
\overline{C}^\bullet _{\Delta }(F_{1[\sigma]})\otimes
 \overline{C}^\bullet _{\Sigma }(F_{2[\tau]})$ we obtain
the required functorial isomorphism
$$\delta :  p_1^{-1}DF_1\otimes p_2^{-1}DF_2 \to DF.$$

b) i) follows from Example 10.3, and ii) follows from the explicit formula for the
morphism $\delta$.
\end{pf}

Consider the diagram of sheaves on $\Phi$
$$\begin{array}{ccccc}
p_1^{-1}\cL _{\Delta} & \otimes & p_2^{-1}\cL _{\Sigma } & = & \cL _{\Phi}\\
p_1^{-1}\epsilon _{\Delta}\downarrow & & \downarrow p_2^{-1}\epsilon _{\Sigma} & &
\downarrow \epsilon _{\Phi}\\
p_1^{-1}D\cL _{\Delta} & \otimes & p_2^{-1}D\cL _{\Sigma } &
\stackrel{\delta}{\rightarrow} & D\cL _{\Phi}.
\end{array}$$
Note that
all arrows are isomorphisms.

\begin{lemma} The above diagram commutes.
\end{lemma}

\begin{pf} The stalks of all the sheaves at the origin $\Or$
are canonically isomorphic to $\bbR$. It follows from Proposition 10.4 b)i) and the definition
of the morphism $\epsilon$ that the diagram commutes at the origin. Hence it commutes by
the rigidity of $\cL _{\Phi}$.
\end{pf}

\begin{prop} Let $\Delta $ and $\Sigma $ be quasi-convex fans, $\Phi =\Delta \times \Sigma $.
 Then the pairing
$$\Gamma (\cL _{\Phi})\times \Gamma _{\Phi ^0}\cL _{\Phi} \to A(2n)$$
is equal to the tensor product of the pairings
$$\Gamma (\cL _{\Delta})\times \Gamma _{\Delta ^0}\cL _{\Delta} \to A_1(2n_1),$$
and
$$\Gamma (\cL _{\Sigma})\times \Gamma _{\Sigma ^0}\cL _{\Sigma} \to A_2(2n_2)$$
under the isomorphisms
$$\Gamma (\cL _{\Phi})=\Gamma (\cL _{\Delta })\otimes
\Gamma (\cL _{\Sigma }),$$
$$\Gamma _{\Phi ^0}\cL _{\Phi}=\Gamma _{\Delta ^0}\cL _{\Delta }\otimes
\Gamma _{\Sigma ^0} \cL _{\Sigma },$$
and
$$A(2n)=A_1(2n_1)\otimes A_2(2n_2).$$
\end{prop}

\begin{pf} Applying the functor of global sections to the commutative
diagram of Lemma 10.5 we obtain the commutative diagram
$$\begin{array}{ccccc}
\Gamma (\cL _{\Delta}) & \otimes & \Gamma (\cL _{\Sigma }) &
\to & \Gamma (\cL _{\Phi})\\
\epsilon _{\Delta}\downarrow & & \downarrow \epsilon _{\Sigma} & &
\downarrow \epsilon _{\Phi}\\
\Hom _{A_1}(\Gamma _{\Delta ^0}\cL _{\Delta},A_1(2n_1)) & \otimes &
\Hom _{A_2}(\Gamma _{\Sigma ^0}\cL _{\Sigma },A_2(2n_2)) &
\stackrel{\delta}{\rightarrow} & \Hom _A(\Gamma _{\Phi ^0}\cL _{\Phi},A(2n))
\end{array}$$
Now apply Proposition 10.4 b) ii).
\end{pf}

\begin{cor} Let $\Delta $ and $\Sigma $ be quasi-convex fans, $\Phi =\Delta \times \Sigma $.
Then the pairing
$$IH(\Phi)\times IH(\Phi ,\partial \Phi)\to \bbR (2n)$$
is the tensor product of pairings
$$IH(\Delta )\times IH(\Delta ,\partial \Delta)\to \bbR (2n_1)$$
and
$$IH(\Sigma)\times IH(\Sigma ,\partial \Sigma)\to \bbR (2n_2).$$
In particular this holds if the fans are complete.
\end{cor}

\subsection{The (HL) and (HR) theorems for product of fans}

Let $l_1$ and $l_2$ be strictly convex piecewise linear functions on the fans $\Delta $ and
$\Sigma $ respectively. Then $l=l_1+l_2$ is strictly convex on $\Phi =\Delta \times \Sigma$.

\begin{thm} Assume that fans $\Delta $, $\Sigma $ (and hence also $\Phi$) are complete.
If the (HL) (resp. (HR)) theorem holds for operators $l_1$ and $l_2$ on $IH(\Delta)$ and
$IH(\Sigma )$ then it also holds for the operator $l$ on $IH(\Phi)$.
\end{thm}

\begin{pf} We have $IH(\Phi)=IH(\Delta )\otimes IH(\Sigma )$ and $l=l_1\otimes 1+
1\otimes l_2$. The statement about the (HL) theorem follows.

Choose a basis for the primitive parts $Prim _{l_1}IH(\Delta)$ and $Prim _{l_2}IH(\Sigma )$
 which is orthogonal with respect to the form $Q_{l_1}$ and $Q_{l_2}$ respectively. This
 defines an orthogonal decomposition of $IH(\Delta)$ and $IH(\Sigma)$ into cyclic $\bbR [l_1]$
 and  $\bbR [l_2]$-modules respectively. But for cyclic modules the assertion can be deduces
 from the classical Hodge-Riemann relations of the primitive cohomology of the product of
 projective spaces.
 (It would be nice to have an elementary algebraic proof but we do not have one.)
\end{pf}

\end{document}